\documentclass[12pt,thmsa]{article}
\usepackage{amsmath}
\usepackage{graphicx}
\usepackage{amsfonts}
\usepackage{amssymb}
\usepackage{mathrsfs}

\newtheorem{lem}{Lemma}[section]

\newtheorem{Rmq}{Remark}[section]

\newtheorem{thm}{Theorem}[section]
\begin{document}

\begin{center}
{\Large \textbf{Kernel based method for the $k$-sample problem}}

\bigskip

Armando Sosthene Kali  BALOGOUN\textsuperscript{a} , Guy Martial  NKIET\textsuperscript{b}  and Carlos OGOUYANDJOU\textsuperscript{a}

\bigskip

\textsuperscript{a}Institut de Math\'ematiques et de Sciences Physiques, Porto Novo, B\'enin.
\textsuperscript{b}Universit\'{e} des Sciences et Techniques de Masuku,  Franceville, Gabon.

\bigskip

E-mail adresses : sosthene.balogoun@imsp-uac.org,  guymartial.nkiet@mathsinfo.univ-masuku.com,  ogouyandjou@imsp-uac.org.

\bigskip
\end{center}

\noindent\textbf{Abstract.}In this paper we deal with the problem of testing for the equality of $k$ probability distributions defined on $(\mathcal{X},\mathcal{B})$, where $\mathcal{X}$ is a metric space and $\mathcal{B}$ is the corresponding Borel $\sigma$-field. We introduce a test statistic based on reproducing kernel Hilbert space embeddings and derive   its asymptotic distribution under the null hypothesis. Simulations show that the introduced procedure outperforms known methods.

\bigskip


\section{Introduction}
\label{Intro}
\noindent Testing for homogeneity, that is  testing for the equality of several probability distributions is an old and important problem in statistics. When the number $k$ of these  distributions is greater than two, it is named the $k$-sample problem and has been  tackled in the literature under different approaches. For instance, the traditional  Kolomogorov-Smirnov, Cram\'er-von Mises  and   Anderson-Darling tests (\cite{con},\cite{gib}), initially introduced to treat the case of two distributions only, have been extended  for dealing with the aforementioned $k$-sample problem  (\cite{kief},\cite{wolf},\cite{sch}).    Also, procedures based the likelihood ratio and which led to more powerful tests  than the previous ones were introduced in \cite{zhang}. Nevertheless,  all these  methods just permit to test the equality of distributions defined on $(\mathbb{R},\mathcal{B}_\mathbb{R})$, where $\mathcal{B}_\mathbb{R}$ is the Borel $\sigma$-field associated to $\mathbb{R}$, and cannot be used for distributions defined on more complex spaces. The interest of  kernel-based  methods, that is methods based on the use of reproducing kernel Hilbert spaces embeddings,  relies on the fact that they permit to deal with high-dimensional and structured data (\cite{harchaoui2}), which the aforementioned traditional  methods do not do. In this vein, Harchaoui et al. \cite{harchaoui1} and, more recently, Gretton et al. \cite{gretton} proposed kernel-based methods for the two sample problem. The former introduced a method based on the  maximum Fisher discriminant ratio  while the latter used   the maximum mean discrepancy.  The extension of their  procedures  to the  case of more than two distributions is of a great interest since, to the best of our knowledge, it it has never been done.

In this paper, we deal with the $k$-sample problem by extending the kernel-based approach of Harchaoui  et al. \cite{harchaoui1}. The rest of the paper is organized as follows.  In Section 2,  we recall some basic  facts about the reproducing kernel Hilbert spaces embeddings.  In Section 3, after specifying the testing problem that we deal with, we introduce a test statistic and derive its asymptotic distribution under the null hypothesis. We also tackle computational aspects that show how to compute this test statistic in practice.   Section 4  is devoted to the presentation of simulations made in order to evaluate performance of our proposal and to compare it with known  methods. All the proofs are postponed in Section 5.

\section{ Preliminary notions}\label{sec2}
\noindent In this section, we recall the notion of reproducing kernel hilbert space (RKHS) and we just define some elements related to it that are useful in this paper. For more details on RKHS and its use in probability and statistics, one may refer to \cite{berlinet}.

\bigskip

\noindent Letting  $(\mathcal{X},\mathcal{B})$ be a measurable space, where $(\mathcal{X},d)$ is a metric space  and $\mathcal{B}$ is the corresponding Borel $\sigma$-field, we consider a Hilbert space $\mathcal{H}$ of functions from $\mathcal{X}$ to $\mathbb{R}$, endowed with an inner product $<\cdot,\cdot>_\mathcal{H}$. This space is said to be a RKHS   if there exists a kernel, that is a symmetric positive semi-definite function $K:\mathcal{X}^2\rightarrow \mathbb{R} $, such that for any $f\in \mathcal{H}$ and any $x\in\mathcal{X}$, one has $K(x,\cdot )\in \mathcal{H}$ and  $f(x)=<f,K(x,\cdot)>_\mathcal{H}$. When $\mathcal{H}$ is a RKHS with kernel $K$,  the map $\Phi\,:\,x\in\mathcal{X}\mapsto K(x,\cdot)\in \mathcal{H}$ characterizes $K$ since one has
\[
K(x,y)=<\Phi(x),\Phi(y)>_{\mathcal{H}}
\]
for any $(x,y)\in\mathcal{X}^2$. It is called the feature map and it is an important tool when dealing with kernel methods for statistical problems.  Throughout this paper, we consider a RKHS $\mathcal{H}$ with kernel $K$ satisfying the following assumptions:

\bigskip

\noindent $(\mathscr{A}_1):$  $\Vert  K\Vert_{\infty}:= \sup\limits_{(x,y)\in \mathcal{X}^2} K(x,y)<+ \infty $;

\bigskip

\noindent $(\mathscr{A}_2):$  the RKHS associated to the kernel $K$ is dense in $L^2(\mathbb{P})$ where $\mathbb{P}$ is a probability measure on $(\mathcal{X},\mathcal{B})$.

\bigskip

\noindent Let   $X$ be a random variable taking values in $\mathcal{X}$ and with probability distribution $\mathbb{P}$. If $\mathbb{E}(\Vert \Phi(X)\Vert_\mathcal{H})=\int_{\mathcal{X}}\Vert\Phi(x)\Vert_\mathcal{H}d\mathbb{P}(x)<+\infty$, the mean element $m$ associated with  $X$ is defined for all functions $f\in \mathcal{H}$ as the unique element in  $\mathcal{H}$ satisfying,
\begin{eqnarray*}
<m,f>_{\mathcal{H}}=\mathbb{E}\left(f(X)\right)=\int_{\mathcal{X}}f(x)d\mathbb{P}(x).
\end{eqnarray*}
Furhermore, if  $\mathbb{E}\left(\Vert \Phi(X)\Vert_\mathcal{H}^2\right)<+\infty$, we can define the covariance operator associated to $X$ as the unique operator $V$ from $\mathcal{H}$ to itself such  that, for any pair $(f,g)\in\mathcal{H}^2$, one has
\begin{eqnarray*}
<f,Vg>_{\mathcal{H}}=\textrm{Cov}\left(f(X),g(X)\right)=\mathbb{E}\left(f(X)g(X)\right)-\mathbb{E}(f(X))\,\mathbb{E}(g(X)).
\end{eqnarray*}
It is very important to note that if   $(\mathscr{A}_1)$ is satisfied, then the mean element $m$  and the covariance operator $V$ are  well-defined. They can also be expressed as
\[
m=\mathbb{E}\left(K(X,\cdot)\right)\,\,\,\,
\]
and
\begin{eqnarray*}
V&=&\mathbb{E}\bigg(\left(K(X,\cdot)-m\right)\otimes \left(K(X,\cdot)-m\right)\bigg)\\
&=&\mathbb{E}\left(K(X,\cdot )\otimes K(X,\cdot)\right)-m\otimes m
\end{eqnarray*}
where $\otimes$ is the tensor product such that, for any pair  $(x,y)\in\mathcal{H}^2$, $x\otimes y$ is the linear map from $\mathcal{H}$ to itself satisfying $(x\otimes y)(h)=<x,h>_\mathcal{H}y$ for all $h\in\mathcal{H}$.
 The   empirical counterparts of $m$ and $V$, obtained from a i.i.d.  sample ${X_1,\cdots,X_n}$    of   $X$,  are  then given by:
\begin{equation*}\label{empirm}
\widehat{m}=\frac{1}{n} \sum_{i=1}^{n}K(X_i,\cdot)
\end{equation*}
and
\begin{equation*}\label{empirv}
\widehat{V}=\frac{1}{n} \sum_{i=1}^{n}\left(K(X_i,\cdot)-\widehat{m}\right)\otimes\left(K(X_i,\cdot)-\widehat{m}\right)
=\frac{1}{n} \sum_{i=1}^{n} K(X_i,\cdot)\otimes K(X_i,\cdot) -\widehat{m}\otimes\widehat{m}.
\end{equation*}

\section{The $k$-sample problem}
\noindent In this section, we specify the $k$-sample problem that we deal with, as a test for hypotheses that are given. Then, a test statistic is proposed and its asymptotic distribution under the null hypothesis is derived. Finally, we deal with computational aspects and show how the introduced test statistic can be computed in practice.

\bigskip

\noindent For $k\in\mathbb{N}^\ast$ such that $k\geq 2$, we consider $k$ probability distibutions $\mathbb{P}_1,\cdots,\mathbb{P}_k$ on $(\mathcal{X},\mathcal{B})$. For $j\in\{1,\cdots,k\}$, we denote by $m_j$  and by $V_j$  the mean element and the covariance operator, respectively, associated to $\mathbb{P}_j$. The $k$-sample problem that we deal with is the test for the hypothesis $\mathscr{H}_0$ : $\mathbb{P}_1=\cdots\mathbb{P}_k$
 against the alternative given by
$\mathscr{H}_1$ : $\exists (j,l)$, $\mathbb{P}_j\neq \mathbb{P}_l$.

\subsection{Test statistic}
\noindent  For $j=1,2,\cdots,k$,  let $\{X_1^{(j)},\cdots,X_{n_j}^{(j)}\}$  be an i.i.d. sample in  $\mathcal{X}$ with commmon distribution  $\mathbb{P}_j$. We consider the statistics 
\[
\widehat{m}_j=\frac{1}{n_j}\sum_{i=1}^{n_j}K(X_i^{(j)},\cdot),
\]
\begin{eqnarray*}
\widehat{V}_j&=&\frac{1}{n_j}\sum_{i=1}^{n_j}\left(K(X_i^{(j)},\cdot)-\widehat{m}_j\right)\otimes\left(K(X_i^{(j)},\cdot)-\widehat{m}_j\right)\\
&=&\frac{1}{n_j}\sum_{i=1}^{n_j}K(X_i^{(j)},\cdot)\otimes K(X_i^{(j)},\cdot)-\widehat{m}_j\otimes \widehat{m}_j
\end{eqnarray*}
from which we define
\[
\widehat{W}_n=\sum_{j=1}^{k}\frac{n_j}{n}\widehat{V}_j,\,\,\,\,\widehat{m}=\sum_{j=1}^{k}\frac{n_j}{n}\,\widehat{m}_j.
\]
where $n=\sum_{j=1}^kn_j$.  Let $\{\gamma_n\}_{n\geq 1}$ be a sequence of strictly positive numbers such that $\lim_{n\rightarrow +\infty}(\gamma_n)=  0$. Then, we consider
\[
\widehat{V}_n=\widehat{W}_n+\gamma_n\mathbb{I},
\]
  where $\mathbb{I}$ denotes the identity operator of $\mathcal{H}$, and we take as test statistic for the $k$-sample problem the statistic:
\begin{eqnarray*}
\widehat{T}_n=\frac{\sum_{j=1}^{k}P_j\parallel\widehat{V}^{-1/2}_n\widehat{\delta}_j\parallel^2_{\mathcal{H}}}{\sqrt{2}\ell(\widehat{W}_n,\gamma_n)},
\end{eqnarray*}
where
\[
P_j=\frac{n_j}{n},\,\,\,\widehat{\delta}_j=\widehat{m}_j-\widehat{m}
\]
and
\begin{eqnarray*}
\ell(\widehat{W}_n,\gamma_n)=\left\{\sum\limits_{p=1}^{+\infty}(\lambda_{p}(\widehat{W}_n)+\gamma_n)^{-2}\lambda^2_{p}(\widehat{W}_n)\right\}^{1/2}.
\end{eqnarray*}
\begin{Rmq}
The quantity $\ell(\widehat{W}_n,\gamma_n)$ is a normalization factor that permits to rescale the statistic  $\sum_{j=1}^{k}P_j\parallel\widehat{V}^{-1/2}_n\widehat{\delta}_j\parallel^2_{\mathcal{H}}$ in order to get a well-grounded test statistic. Note that in \cite{harchaoui1} a factor for recentering is also introduced, but we do not need it in this paper. It is know from \cite{harchaoui2}  that  
\begin{equation}\label{ell}
\ell(\widehat{W}_n,\gamma_n)=\textrm{tr}\left(\widehat{V}_n^{-2}\widehat{W}_n^2\right).
\end{equation}
\end{Rmq}

\subsection{Asymptotic distribution under $\mathscr{H}_0$}
  We consider the following assumptions:
\bigskip

\noindent $(\mathscr{A}_3):$  For $j\in\{1,\cdots,k\}$, one has $\lim_{n_j\rightarrow +\infty}\frac{n_j}{n}=\rho_j$, where $\rho_j$ is a real belonging to $]0,1[$.

\bigskip

\noindent $(\mathscr{A}_4):$  the eigenvalues $\{\lambda_p(V_j)\}_{p\geq 1}$ satisfy $\sum_{p=1}^{+\infty}\lambda_p^{1/2}(V_j)<+\infty$ for $j=1,2,\cdots,k$;

\bigskip

\noindent $(\mathscr{A}_5):$ there are infinitely many strictly positive eigenvalues $\{\lambda_p(V_j)\}_{p\geq1 }$ of $V_j$ for $j=1,2,\cdots,k$. 

\bigskip

\noindent Then, we have:

\begin{thm}\label{loi}
	Assume ($\mathscr{A}_1$)  to ($\mathscr{A}_5$) and that $\lim_{n\rightarrow +\infty} (\gamma_n+\gamma^{-1}_nn^{-1/2})= 0$, then  under $\mathscr{H}_0$,  $n\widehat{T}_n$ converges in distribution, as $n\rightarrow +\infty$,  to  $\mathcal{N}(0,1)$.
	
\end{thm}
\subsection{Computation of the test statistic}

\noindent For computing this test statistic in practice, the kernel trick (\cite{shawe}) can be used as it was already done in \cite{harchaoui1} for twe two-groups case. For $j=1,\cdots,k$, we consider the operator $G_n^{(j)}$ from $\mathbb{R}^{n_j}$ to $\mathcal{H}$ represented in matrix form as
\begin{align*}
G^{(j)}_{n}=[K(X^{(j)}_{1},.),\cdots,K(X^{(j)}_{n_{j}},.)].
\end{align*}
Then put $
G_{n}=[G^{(1)}_{n}G^{(2)}_{n}\cdots G^{(k)}_n]$, and consider
\[
 \Lambda_n^{(j,l)} =(G_n^{(j)})^TG_n^{(l)} =\left(
\begin{array}{cccc}
K(X_1^{(j)},X_1^{(l)}) & K(X_1^{(j)},X_2^{(l)}) & \cdots & K(X_1^{(j)},X_{n_l}^{(l)}) \\
K(X_2^{(j)},X_1^{(l)}) & K(X_2^{(j)},X_2^{(l)}) & \cdots &K(X_2^{(j)},X_{n_l}^{(l)}) \\
\vdots & \vdots & \ddots & \vdots \\
K(X_{n_j}^{(j)},X_1^{(l)}) & K(X_{n_j}^{(j)},X_2^{(l)}) & \cdots &K(X_{n_j}^{(j)},X_{n_l}^{(l)}) \\
\end{array}
\right)
\]
and the Gram matrix
\[
\Lambda_n=G_n^TG_n=\left(
\begin{array}{cccc}
\Lambda_n^{(1,1)} & \Lambda_n^{(1,2)} & \cdots & \Lambda_n^{(1,k)} \\
\Lambda_n^{(2,1)} & \Lambda_n^{(2,2)} & \cdots & \Lambda_n^{(2,k)} \\
\vdots & \vdots & \ddots & \vdots \\
\Lambda_n^{(k,1)} & \Lambda_n^{(k,2)} & \cdots & \Lambda_n^{(k,k)} \\
\end{array}
\right).
\]
Further, denoting by $\textrm{\textbf{I}}_l$ (resp.   $\textrm{\textbf{1}}_l$) the $l\times l$ identity matrix (resp. the $l\times 1$ vector whose components are all equal to $1$), we consider the matrices $\textrm{\textbf{Q}}_{n_j} =\textrm{\textbf{I}}_{n_j}-n_j^{-1}\textrm{\textbf{1}}_{n_j}\textrm{\textbf{1}}_{n_j}^T$,
\begin{eqnarray*}\label{1}
N_n = \left(
\begin{array}{cccc}
\textrm{\textbf{Q}}_{n_1} & 0 & \cdots & 0 \\
0 &\textrm{\textbf{Q}}_{n_2} & \cdots &0 \\
\vdots & \vdots & \ddots & \vdots \\
0 & 0 & \cdots &\textrm{\textbf{Q}}_{n_k} \\
\end{array}
\right),
\end{eqnarray*}
and the vector 
\[
m_{n}^{(j)}=\left(
\begin{array}{c}
m_{n,1}^{(j)} \\
\vdots  \\
m_{n,n}^{(j)} \\
\end{array}
\right)
\]
such that
\[
m_{n,i}^{(j)}=
\left\{
\begin{array}{lcl}
-n^{-1}  &  &\textrm{ if  }\,\,\,  1 \leq i \leq \nu_{j-1}\\
n_j^{-1}-n^{-1} & &\textrm{ if  }\,\,\, \nu_{j-1} +1\leq i \leq  \nu_{j}\\
-n^{-1}& &\textrm{ if  }\,\,\,  \nu_{j}+1\leq i \leq n  
\end{array}\right.
\]
where $\nu_l=\sum_{s=1}^{l}n_l$. Clearly,
\begin{eqnarray*}
\widehat{\delta}_j&=&\widehat{m}_j-\widehat{m}\nonumber=\widehat{m}_j-\sum_{l=1}^{k}P_l\widehat{m}_l\nonumber\\&=&-P_1\widehat{m}_1-P_2\widehat{m}_2-\cdots-P_{j-1}\widehat{m}_{j-1}+(n_j^{-1}-n^{-1})\widehat{m}_j-P_{j+1}\widehat{m}_{j+1}-\cdots-P_k\widehat{m}_k\nonumber\\&=&-n^{-1}\sum_{i=1}^{n_1}K(X_i^{(1)},.)-n^{-1}\sum_{i=1}^{n_2}K(X_i^{(2)},.)-\cdots-n^{-1}\sum_{i=1}^{n_{j-1}}K(X_i^{(j-1)},.)\nonumber\\&+&((n_j^{-1}-n^{-1})\sum_{i=1}^{n_j}K(X_i^{(j)},.)-n^{-1}\sum_{i=1}^{n_{j+1}}K(X_i^{(j+1)},.)-\cdots-n^{-1}\sum_{i=1}^{n_k}K(X_i^{(k)},.)\nonumber\\&=& G_{n}m_{n}^{(j)}
\end{eqnarray*}
and, as in \cite{harchaoui1}, 
$\widehat{V}_{j}=n_j^{-1}G_n^{(j)}\textrm{\textbf{Q}}_{n_j}\textrm{\textbf{Q}}_{n_j}^T(G_n^{(j)})^T$. Therefore 
\begin{eqnarray}\label{wn}
\widehat{W}_n&=&n^{-1}\sum_{j=1}^{k}G_n^{(j)}\textrm{\textbf{Q}}_{n_j}\textrm{\textbf{Q}}_{n_j}^T(G_n^{(j)})^T=n^{-1}G_nN_nN_n^TG_n^T
\end{eqnarray}
and $
\parallel\widehat{V}^{-1/2}_n\widehat{\delta}_j\parallel^2_{\mathcal{H}}=(m_n^{(j)})^TG_n^T\left(\gamma_n\mathbb{I}+n^{-1}G_nN_nN_n^TG_n^T\right)^{-1}G_nm_n^{(j)}$. Using the matrix inversion lemma, as in \cite{harchaoui1}, we obtain
\begin{eqnarray}
&&P_j\parallel\widehat{V}^{-1/2}_n\widehat{\delta}_j\parallel^2_{\mathcal{H}}\nonumber\\
&=&P_j(m_n^{(j)})^TG_n^T\gamma_n^{-1}\{\textrm{\textbf{I}}_{n}-G_nN_n(\gamma_n\textrm{\textbf{I}}_{n}+n^{-1}N_n^TG_n^TG_nN_n)^{-1}N_n^TG_n^T\}G_nm_n^{(1)}\nonumber\\
&=&\frac{n_j}{n\gamma_n}\bigg\{(m_n^{(j)})^T\Lambda_nm_n^{(j)}-n^{-1}(m_n^{(j)})^T\Lambda_nN_n(\gamma_n\textrm{\textbf{I}}_{n}+n^{-1}N_n^T\Lambda_nN_n)^{-1}N_n^T\Lambda_nm_n^{(j)}\bigg\}\nonumber.
\end{eqnarray}
Finally,
\begin{eqnarray*}
\sum_{j=1}^{k}P_j\parallel\widehat{V}^{-1/2}_n\widehat{\delta}_j\parallel^2_{\mathcal{H}}&=&\sum_{j=1}^{k}\frac{n_j}{n\gamma_n}\bigg\{(m_n^{(j)})^T\Lambda_nm_n^{(j)}\nonumber\\
& &-n^{-1}(m_n^{(j)})^T\Lambda_nN_n(\gamma_n\textrm{\textbf{I}}_{n}+n^{-1}N_n^T\Lambda_nN_n)^{-1}N_n^T\Lambda_nm_n^{(j)}\bigg\}.
\end{eqnarray*} 
For computing $\ell(\widehat{W}_n,\gamma_n)$ we use (\ref{ell}) and (\ref{wn}); putting $H_n=n^{-1}G_nN_nN_n^TG_n^T$ we have
$
\ell(\widehat{W}_n,\gamma_n)=\textrm{tr}\left(\left(H_n+\gamma_n\mathbb{I}\right)^{-2}H_n^2\right)$. 
Clearly, $H_n^2=n^{-2}G_nN_n^T\Lambda_nN_nN_n^TG_n^T$ and using the matrix inversion lemma, we obtain
\[
\left(H_n+\gamma_n\mathbb{I}\right)^{-1}=\gamma_n^{-1}\mathbb{I}-n^{-1}\gamma_n^{-1}G_nN_nM_nN_n^TG_n^T
\]
where $M_n=\left(\textrm{\textbf{I}}_{n}+n^{-1}\gamma^{-1}_nN_n^T\Lambda_n N_n\right)^{-1}$. Hence
\begin{eqnarray*}
\left(H_n+\gamma_n\mathbb{I}\right)^{-2}H_n^2&=&\gamma_n^{-2}H_n^2-2n^{-1}\gamma_n^{-3}G_nN_nM_nN_n^TG_n^TH_n^2\\
& &+n^{-2}\gamma_n^{-4}G_nN_nM_nN_n^T\Lambda_n N_nM_nN_n^TG_n^TH_n^2\\
&=&\gamma_n^{-2}H_n^2-2n^{-3}\gamma_n^{-3}G_nN_nM_n\left(N_n^T\Lambda_nN_n\right)^2N_n^TG_n^T\\
& &+n^{-4}\gamma_n^{-4}G_nN_n(M_nN_n^T\Lambda_n N_n)^2N_n^T\Lambda_nN_nN_n^TG_n^T\\
\end{eqnarray*}
and using the property $\textrm{tr}(AB)=\textrm{tr}(BA)$, we finally obtain
\begin{eqnarray*}
\ell(\widehat{W}_n,\gamma_n)&=&\textrm{tr}\bigg(n^{-2}\gamma^{-2}_n(N_n^T\Lambda_n N_n)^2-2n^{-3}\gamma_n^{-3}M_n(N_n^T\Lambda_n N_n)^3\\
& &\hspace{2cm}+n^{-4}\gamma_n^{-4}(M_nN_n^T\Lambda_n N_n)^2(N_n^T\Lambda_n N_n)^2\bigg).
\end{eqnarray*}
\section{Power comparison by Monte Carlo simulation}
\noindent In this section, the empirical power  of the proposed test is computed  through   Monte Carlo simulations and compared to that of tests introduced by Zhang and Wu \cite{zhang} which are based on statistics denoted by $Z_a$,  $Z_c$ and  $Z_k$ obtained from the likelihood-ratio test statistic  and shown to be   more powerful than the classical Kolmogorov-Smirnov, Cram\'er-von Mises and Anderson-Darling $k$-sample tests. We  estimate the powers of  our test and the three  aforementioned  tests   in the following cases ($k = 3$): 
\begin{eqnarray*}
&&\textrm{\textbf{Case 1:} } \mathbb{P}_1=N(3,1),\hspace*{0.2cm}\mathbb{P}_2=Gamma(3,1)\hspace*{0.2cm}\textrm{ and }\hspace*{0.2cm} \mathbb{P}_3=Gamma(6,2);\\
&&\textrm{\textbf{Case 2:} } \mathbb{P}_1=N(0,1),\hspace*{0.2cm}\mathbb{P}_2=N(0,2)\hspace*{0.2cm}\textrm{ and }\hspace*{0.2cm} \mathbb{P}_3=N(0,4);\\
&&\textrm{\textbf{Case 3:} } \mathbb{P}_1=Uniform(0,1),\hspace*{0.2cm}\mathbb{P}_2=Beta(1,1.5)\hspace*{0.2cm}\textrm{ and } \hspace*{0.2cm} \mathbb{P}_3=Beta(1.5,1);\\
&&\textrm{\textbf{Case 4:} } \mathbb{P}_1=N(0,1),\hspace*{0.2cm}\mathbb{P}_2=N(0.3,1)\hspace*{0.2cm}\textrm{ and }\hspace*{0.2cm} \mathbb{P}_3=N(0.6,1).
\end{eqnarray*}
For all tests we take  the significance level $\alpha= 0.05$ and the empirical power is computed over $100$ independent replications. For our test, we used the gaussian  kernel $K(x,y)=\exp[-2(x-y)^2]$, and computed the test statistic as indicated in Section 3.3 by taking 
\[
\gamma_n=0.2\times\textrm{\textbf{1}}_{[1,100 [}(n)+0.01\times\textrm{\textbf{1}}_{[100,300 ]}(n)+n^{-0.25}\times\textrm{\textbf{1}}_{]300,+\infty ]}(n).
\]
The results are given in Figures 1 to 4  that   plot  the  empirical power versus the total sample size $n=n_1+n_2+n_3$. They show that our test outperforms the three tests of Zhang and Wu \cite{zhang} in all cases.

\begin{figure}[h]
	\begin{minipage}[h]{.45\linewidth}
		\centering\includegraphics[width = \linewidth]{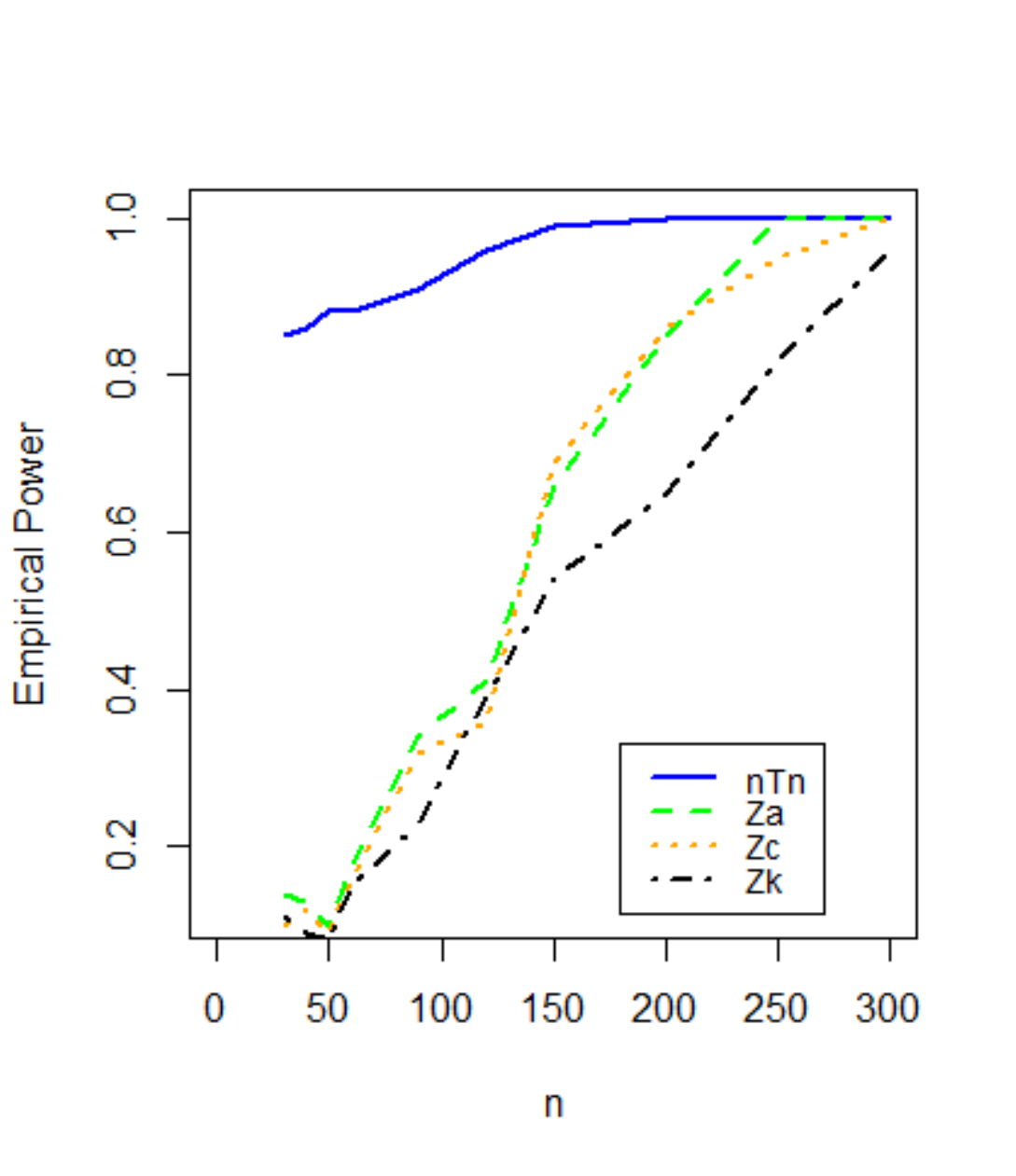}
\caption{Empirical power    versus $n$ for Case 1 with significance level $\alpha=0.05$}
	\end{minipage}
	\hfill
	\begin{minipage}[h]{.45\linewidth}
		\centering\includegraphics[width = \linewidth]{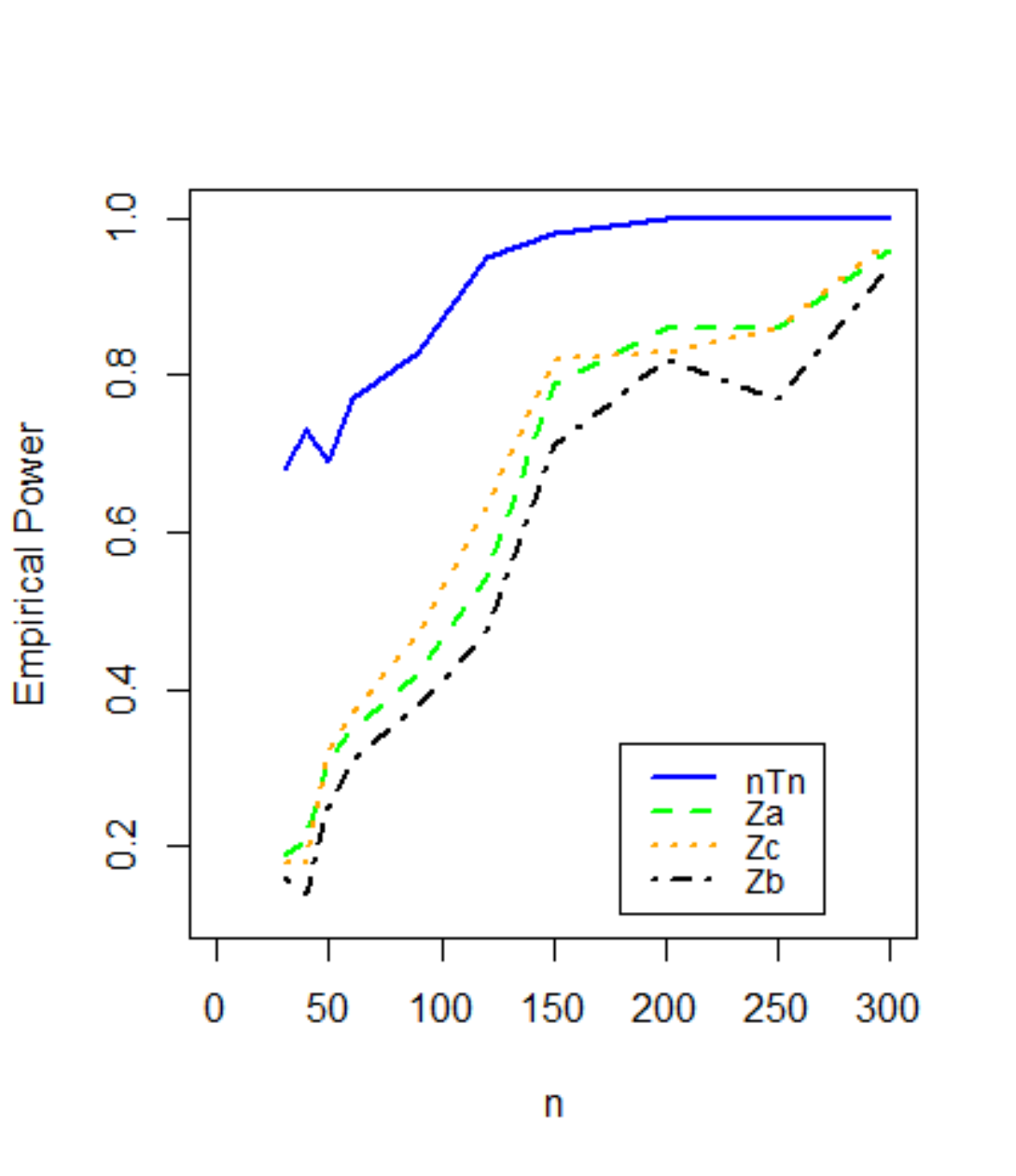}
\caption{Empirical power    versus $n$ for Case 2 with significance level $\alpha=0.05$}
	\end{minipage}
	\label{etiquette_figC}
\end{figure}
\begin{figure}[h]
	\begin{minipage}[h]{.45\linewidth}
		\centering\includegraphics[width = \linewidth]{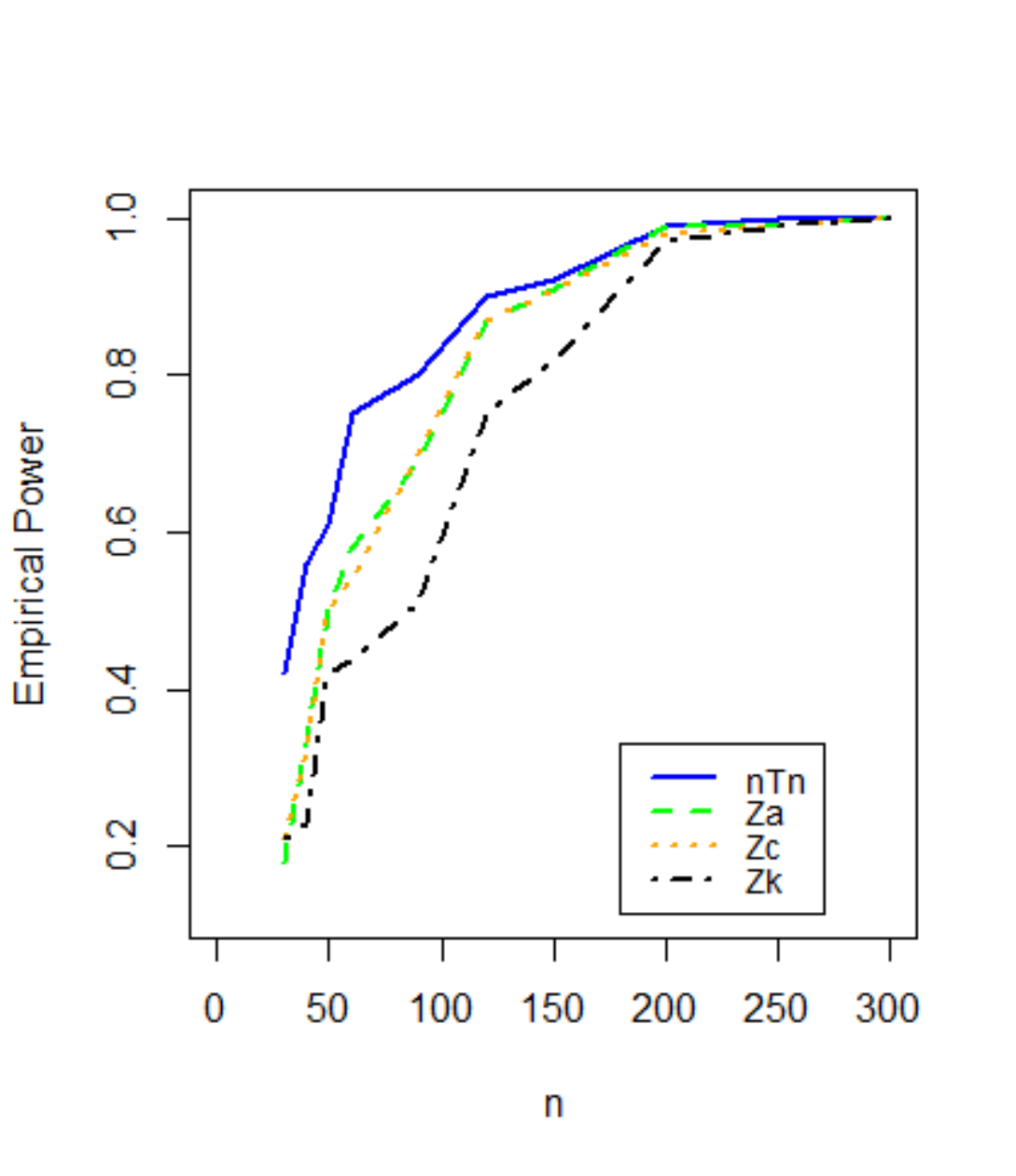}
\caption{Empirical power    versus $n$ for Case 3  with significance level $\alpha=0.05$}
	\end{minipage}
	\hfill
	\begin{minipage}[h]{.45\linewidth}
		\centering\includegraphics[width = \linewidth]{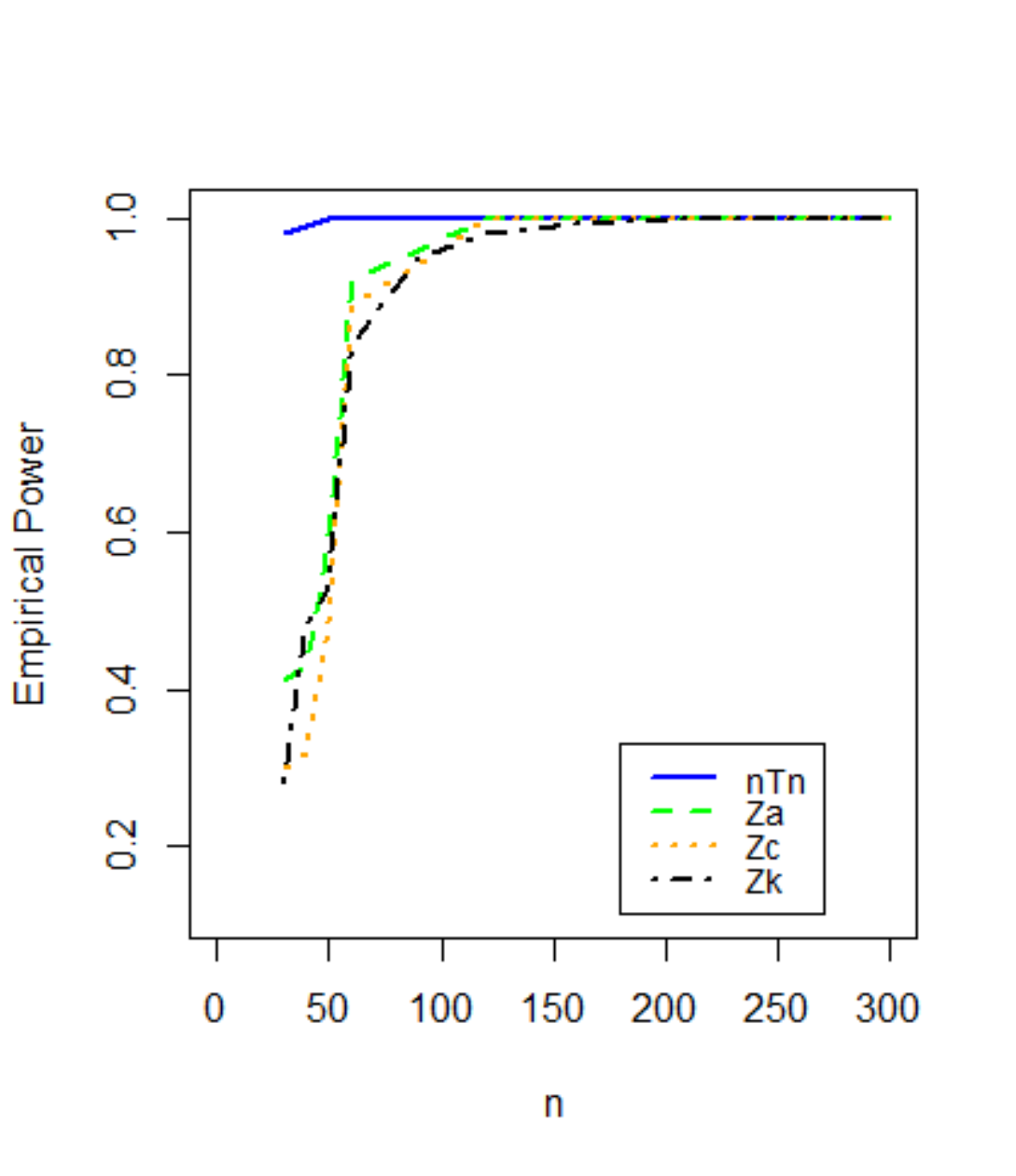}
\caption{Empirical power    versus $n$ for Case 4 with significance level $\alpha=0.05$}
	\end{minipage}
	\label{etiquette_figC}
\end{figure}

\section{Proofs}
\subsection{Preliminary results}
\noindent In this section, we give some results that are necessary for proving Theorem \ref{loi}. 
\begin{lem}\label{l1}
	Assume ($\mathscr{A}_1$), ($\mathscr{A}_3$) and ($\mathscr{A}_4$) . Then, putting $W=\sum_{j=1}^kP_jV_j$, we have $
	 \parallel \widehat{W}_n-W\parallel_{HS}=O_P(n^{-1/2})$.
\end{lem}
\noindent\textit{Proof}. Let   $\{e_p\}_{p\geq 1}$ be an orthonormal basis of $\mathcal{H}$ consisting of eigenvectors of  $V_j$ such that $e_p$ is associated  to the $p$-th  eigenvalue $\lambda_p(V_j)$. Using  Lemma 21 in \cite{harchaoui1} and the equality  $\widehat{W}_n-W=\sum_{j=1}^{k}\frac{n_j}{n}(\widehat{V}_j-V_j)$, we obtain

	\begin{eqnarray*}\label{ineq1}
	\sum_{p=1}^{+\infty}\mid \lambda_p(\widehat{W}_n-W)\mid\leq \sum_{p=1}^{+\infty}\parallel(\widehat{W}_n-W)e_p\parallel_{\mathcal{H}}\leq \sum_{j=1}^{k}\frac{n_j}{n}\left\{\sum_{p=1}^{+\infty}\parallel(\widehat{V}_j-V_j)e_p\parallel_{\mathcal{H}}\right\}.
	\end{eqnarray*} 
	From  Proposition 12 of \cite{harchaoui1}, we have $\sum_{p=1}^{+\infty}\parallel(\widehat{V}_j-V_j)e_p\parallel_{\mathcal{H}}=O_P(n^{-1/2})$ for any $j\in\{1,\cdots,k\}$,
	and since $\lim_{n_j\rightarrow +\infty }(\frac{n_j}{n})=\rho_j$  it follows from  (\ref{ineq1}) that we have the equality $
\sum_{p=1}^{+\infty}\mid \lambda_p(\widehat{W}_n-W)\mid=O_P(n^{-1/2})$.
	Furthermore,
	
\begin{eqnarray*}
\parallel \widehat{W}_n-W\parallel_{HS}&=&\left[\sum\limits_{p=1}^{+\infty}\parallel(\widehat{W}_n-W)e_{p}\parallel_{\mathcal{H}}^2\right]^{1/2}\nonumber\\&\leq&	\sum\limits_{p=1}^{+\infty}\parallel(\widehat{W}_n-W)e_{p}\parallel_{\mathcal{H}}\nonumber\\&\leq&\sum_{j=1}^{k}\frac{n_j}{n}\left\{\sum_{p=1}^{+\infty}\parallel(\widehat{V}_j-V_j)e_p\parallel_{\mathcal{H}}\right\}
\end{eqnarray*}
which proves that   $\parallel\widehat{W}_n-W\parallel_{HS}=O_{P}(n^{-1/2})$.
\hfill $\Box$

\bigskip
\noindent The following lemma  gives  an asymptotic approximation of the test statistic.
\begin{lem}\label{l2}
	Assume ($\mathscr{A}_1$), ($\mathscr{A}_3$) and ($\mathscr{A}_4$). If $\lim_{n\rightarrow +\infty}\left(\gamma_n+\gamma^{-1}_nn^{-1/2}\right)=0$, then
	\begin{eqnarray}\label{resl2}
	n\widehat{T}_n&=&\widetilde{S}_n+o_P(1),
\end{eqnarray}
	 where
\[
	 \widetilde{S}_n=\frac{\sum_{j=1}^{k}n_j\parallel V^{-1/2}_n\widehat{\delta}_j\parallel^2_{\mathcal{H}}}{\sqrt{2}\ell(W,\gamma_n)}.
	\]
	\end{lem} 
   
\noindent\textit{Proof.}
Using Lemma 23 in \cite{harchaoui1}, we have
	\[
\mid\ell(\widehat{W}_n,\gamma_n)-\ell(W,\gamma_n)\mid\leq \frac{\gamma_n^{-1}\parallel \widehat{W}_n-W\parallel_{HS}}{1-\gamma_n^{-1}\parallel \widehat{W}_n-W\parallel_{HS}}. 
\]
Then, from Lemma \ref{l1}  it follows
	\begin{equation}\label{op1}
\mid\ell(\widehat{W},\gamma_n)-\ell(W,\gamma_n)\mid=O_P(\gamma^{-1}_nn^{-1/2})=o_P(1).
\end{equation} 
Furthermore,
		\begin{eqnarray}\label{ineq11}
			\left| \parallel\widehat{V}_n^{-1/2}\widehat{\delta}_j\parallel_{\mathcal{H}}^2-\parallel V_n^{-1/2}\widehat{\delta}_j\parallel_{\mathcal{H}}^2\right|
 &=&\mid <\widehat{\delta}_j,\widehat{V}_n^{-1}\widehat{\delta}_j>_{\mathcal{H}}-<\widehat{\delta}_j,V_n^{-1}\widehat{\delta}_j>_{\mathcal{H}}\mid \nonumber\\
&=&\mid <\widehat{V}_n^{-1}\widehat{\delta}_j,(V_n-\widehat{V}_n)V_n^{-1}\widehat{\delta}_j>_{\mathcal{H}}\mid \nonumber\\
&\leq&\parallel \widehat{V}_n^{-1}\widehat{\delta}_j \parallel_{\mathcal{H}} \,\parallel V_n-\widehat{V}_n\parallel_{HS} \, \parallel V_n^{-1}\widehat{\delta}_j\parallel_{\mathcal{H}},
		\end{eqnarray}
	
\begin{eqnarray}\label{ina1}
	\parallel V_n^{-1}\widehat{\delta}_j\parallel_{\mathcal{H}}^2&=&<\widehat{\delta}_j,V_n^{-2}\widehat{\delta}_j>_{\mathcal{H}}\nonumber=\textrm{tr}\left(\widehat{\delta}_j\otimes (V_n^{-2} \widehat{\delta}_j)\right)=\textrm{tr}\left(V_n^{-2}(\widehat{\delta}_j\otimes \widehat{\delta}_j)\right) \nonumber\\&=&\sum_{p=1}^{+\infty}<V_n^{-2}e_p,(\widehat{\delta}_j\otimes\widehat{\delta}_j)e_p>_{\mathcal{H}}\nonumber\\&=&\sum_{p=1}^{+\infty}(\lambda_p+\gamma_n)^{-2}<e_p,\widehat{\delta}_j>_{\mathcal{H}}^2\nonumber\\&\leq& \sum_{p=1}^{+\infty}(\lambda_p+\gamma_n)^{-2}\parallel e_p \parallel_{\mathcal{H}}^2 \parallel \widehat{\delta}_j\parallel_{\mathcal{H}}^2 \nonumber\\&\leq& \sum_{p=1}^{+\infty}(\lambda_p+\gamma_n)^{-2} \parallel \widehat{\delta}_j\parallel_{\mathcal{H}}^2. 
\end{eqnarray}
and
\begin{eqnarray}\label{indel}
\parallel \widehat{\delta}_j\parallel_{\mathcal{H}}&=&\parallel \sum_{\underset{l\neq j}{l=1}}^{k}P_l(\widehat{m}_j-\widehat{m}_l)\parallel_{\mathcal{H}}\leq\sum_{\underset{l\neq j}{l=1}}^{k}P_l\parallel \widehat{m}_j-\widehat{m}_l\parallel_{\mathcal{H}}\nonumber\\& \leq&\sum_{\underset{l\neq j}{l=1}}^{k}P_l\left(\parallel \widehat{m}_j-m_j\parallel_{\mathcal{H}} +\parallel \widehat{m}_l-m_l\parallel_{\mathcal{H}}\right).
\end{eqnarray}
 Using the central limit theorem, we have $\parallel \widehat{m}_l-m_l\parallel_{\mathcal{H}}=O_P(n^{-1/2})$ for all  $l\in \{1,2,\cdots,k\}$, and since $\lim\limits_{n_l\rightarrow +\infty}(\frac{n_l}{n})=\rho_l$ it follows from (\ref{indel}) that   $\parallel \widehat{\delta}_j\parallel_{\mathcal{H}}^2=o_P(1)$.  The fact that
   $\sum_{p=1}^{+\infty}(\lambda_p+\gamma_n)^{-2}<+\infty$ permits to deduce from (\ref{ina1}) that   $\parallel V_n^{-1}\widehat{\delta}_j\parallel_{\mathcal{H}}^2=o_p(1).$	
Moreover,
	 \begin{eqnarray}\label{ina2}
	 	\parallel \widehat{V}_n^{-1}\widehat{\delta}_j \parallel_{\mathcal{H}}&\leq& \parallel{V}_n^{-1}\widehat{\delta}_j\parallel_{\mathcal{H}}+\parallel (\widehat{V}_n^{-1}- {V}_n^{-1})\widehat{\delta}_j \parallel_{\mathcal{H}}\nonumber\\& =& \parallel{V}_n^{-1}\widehat{\delta}_j\parallel_{\mathcal{H}}+\parallel \widehat{V}_n^{-1}(V_n-\widehat{V}_n ){V}_n^{-1}\widehat{\delta}_j \parallel_{\mathcal{H}}\nonumber\\&\leq& \parallel{V}_n^{-1}\widehat{\delta}_j\parallel_{\mathcal{H}}+\parallel \widehat{V}_n^{-1}\parallel_{HS}\parallel V_n-\widehat{V}_n \parallel_{HS} \parallel{V}_n^{-1}\widehat{\delta}_j\parallel_{\mathcal{H}}
	 	\end{eqnarray}

	 \noindent Next, using the upper-bound $\left\| \widehat{V}^{-1}_n\right\|_{HS}\leq \gamma^{-1}_n$ and  Lemma \ref{l1},
\begin{equation}\label{ineq3}
\parallel \widehat{V}_n-V\parallel_{HS}=\parallel \widehat{W}_n-W\parallel_{HS}=O_P(n^{-1/2})
\end{equation}
we get
	 $\parallel \widehat{V}_n^{-1}\parallel_{HS}\parallel V_n-\widehat{V}_n \parallel_{HS}=O_p(\gamma^{-1}_nn^{-1/2})=o_P(1)$. Since we have the equality   $\parallel{V}_n^{-1}\widehat{\delta}_j\parallel_{\mathcal{H}}=o_P(1)$ it follows from (\ref{ina2}) that  $\parallel \widehat{V}_n^{-1}\widehat{\delta}_j \parallel_{\mathcal{H}}=o_P(1)$. Finally,  using all that precedes together with (\ref{ineq11}) and (\ref{ineq3}), we conclude that
\begin{equation}\label{op2}
 \left| \parallel\widehat{V}_n^{-1/2}\widehat{\delta}_j\parallel_{\mathcal{H}}^2-\parallel V_n^{-1/2}\widehat{\delta}_j\parallel_{\mathcal{H}}^2\right|=o_P(1).
\end{equation}
Now, we will use the above results for proving  (\ref{resl2}). We have
\begin{eqnarray}\label{op3}
\left \vert n\widehat {T}_n-\widetilde{S}_n\right\vert&\leq&\left\vert\frac{1}{\ell(\widehat{W}_n,\gamma_n)}-\frac{1}{\ell(W,\gamma_n)}\right\vert
\sum_{j=1}^{k}n_j\parallel\widehat{V}^{-1/2}_n\widehat{\delta}_j\parallel^2_{\mathcal{H}}\nonumber\\
& &+\frac{1}{\ell(W,\gamma_n)}\sum_{j=1}^{k}n_j\left\vert\parallel V^{-1/2}_n\widehat{\delta}_j\parallel^2_{\mathcal{H}}- \parallel\widehat{V}^{-1/2}_n\widehat{\delta}_j\parallel^2_{\mathcal{H}}\right\vert\nonumber\\
&\leq &\frac{\left\vert\ell(\widehat{W}_n,\gamma_n)-\ell(W,\gamma_n)\right\vert}{\ell(\widehat{W}_n,\gamma_n)\ell(W,\gamma_n)}
\sum_{j=1}^{k}n_j\left\vert\parallel \widehat{V}^{-1/2}_n\widehat{\delta}_j\parallel^2_{\mathcal{H}}- \parallel V^{-1/2}_n\widehat{\delta}_j\parallel^2_{\mathcal{H}}\right\vert\nonumber\\
& &+\frac{\left\vert\ell(\widehat{W}_n,\gamma_n)-\ell(W,\gamma_n)\right\vert}{\ell(\widehat{W}_n,\gamma_n)\ell(W,\gamma_n)}
\sum_{j=1}^{k}n_j \parallel V^{-1/2}_n\widehat{\delta}_j\parallel^2_{\mathcal{H}}\nonumber\\
& &+\frac{1}{\ell(W,\gamma_n)}\sum_{j=1}^{k}n_j\left\vert\parallel \widehat{V}^{-1/2}_n\widehat{\delta}_j\parallel^2_{\mathcal{H}}- \parallel V^{-1/2}_n\widehat{\delta}_j\parallel^2_{\mathcal{H}}\right\vert
\end{eqnarray}
and since $ \parallel V^{-1/2}_n\widehat{\delta}_j\parallel^2_{\mathcal{H}}\leq \parallel \widehat{\delta}_j\parallel_{\mathcal{H}}\,\,\parallel V_n^{-1}\widehat{\delta}_j\parallel_{\mathcal{H}}$, it follows that $ \parallel V^{-1/2}_n\widehat{\delta}_j\parallel^2_{\mathcal{H}}=o_P(1)$. Therefore, from (\ref{op1}), (\ref{op2}) and (\ref{op3}) we deduce that $\left \vert n\widehat {T}_n-\widetilde{S}_n\right\vert$ converges in probabilty to $0$ as $n\rightarrow +\infty$.
\hfill $\Box$

\bigskip

\noindent For ease of notation, in the following, we shall denote by  $\lambda_p$, $\lambda_r$, $\ell_n$ the terms  $\lambda_p(W)$, $\lambda_r(W)$, $\ell(W,\gamma_n)$ respectively. Define
\begin{eqnarray*}\label{a1}
 Y_{n,p,i,j,l}=\left\{\begin{array}{lc}
  \frac{P_k\sqrt{n_j}}{n_{j}}\left(e_{p}(X^{(j)}_{i})-\mathbb{E}[e_{p}(X^{(j)}_{1})]\right) \,\,\,\,\,\,1\leq i \leq n_{j}\\
-\frac{P_k\sqrt{n_j}}{n_{l}}\left(e_{p}(X^{(l)}_{i-n_{j}})-\mathbb{E}[e_{p}(X^{(l)}_{1})]\right) \,\,\,\,\,\,n_{j}+1\leq i \leq n_j+n_l.  
\end{array}\right.
\end{eqnarray*}
Then, we have:
\begin{lem}\label{l3}We have:
	\begin{eqnarray}
	&&\sum_{j=1}^{k}\sum_{\stackrel{l=1}{l\neq j}}^{k}\sum_{i=1}^{n_j+n_l}\mathbb{E}\left[Y_{n,p,i,j,l}Y_{n,r,i,j,l}\right]=\lambda_p^{1/2}\lambda_r^{1/2}\left[1+O(n^{-1})\right]\delta_{pr}\label{d};\\
	&&Cov\left(Y_{n,p,i,j,l}^2,Y_{n,r,i,j,l}^2\right)\leq Cn^{-2}\Vert K \Vert _{\infty}\lambda_p^{1/2}\lambda_r^{1/2}\label{f},
	\end{eqnarray}   
where $C$ is a positive constant.
\end{lem}
\noindent\textit{Proof.}
	From elementary calculation we have
	\begin{eqnarray}
	&&\sum_{i=1}^{n_{j}}\mathbb{E}[Y_{n,p,i,j,l}Y_{n,r,i,j,l}]=\frac{n_{l}^2}{n^2}\lambda^{1/2}_{p}\lambda^{1/2}_{r}\delta_{pr},\nonumber\\
&&\sum_{i=n_{j}+1}^{n_j+n_l}\mathbb{E}[Y_{n,p,i,j,l}Y_{n,r,i,j,l}]=\frac{n_{j}n_l}{n^2}\lambda^{1/2}_{p}\lambda^{1/2}_{r}\delta_{pr}.\nonumber
	\end{eqnarray}
	Hence 
	\begin{eqnarray*}
\sum_{j=1}^{k}\sum_{\stackrel{l=1}{l\neq j}}^{k}\sum_{i=1}^{n_j+n_l}\mathbb{E}\left[Y_{n,p,i,j,l}Y_{n,r,i,j,l}\right]&=&\sum_{j=1}^{k}\sum_{\stackrel{l=1}{l\neq j}}^{k}\frac{n_l(n_l+n_j)}{n^2}\lambda_p^{1/2}\lambda_r^{1/2}\delta_{pr}\nonumber\\&=&\sum_{j=1}^{k}\sum_{\stackrel{l=1}{l\neq j}}^{k}\left[P_l^2+P_jP_l\right]\lambda_p^{1/2}\lambda_r^{1/2}\delta_{pr}\nonumber\\&=&\left[1+k\sum_{l=1}^{k}P_l^2-2\sum_{l=1}^{k}P_l^2\right]\lambda_p^{1/2}\lambda_r^{1/2}\delta_{pr}\nonumber\\&=&\left[1+(k-2)\sum_{l=1}^{k}P_l^2\right]\lambda_p^{1/2}\lambda_r^{1/2}\delta_{pr}\nonumber\\&=&\left[1+\frac{k-2}{n^2}\sum_{l=1}^{k}n_l^2\right]\lambda_p^{1/2}\lambda_r^{1/2}\delta_{pr}\nonumber\\&=&\left[1-\frac{k-2}{n}+\frac{2k-4}{n^2}\sum_{l=1}^{k}\sum_{i=1}^{n_l}i\right]\lambda_p^{1/2}\lambda_r^{1/2}\delta_{pr}\nonumber\\&=&\left[1+O(n^{-1})\right]\lambda_p^{1/2}\lambda_r^{1/2}\delta_{pr}.
	\end{eqnarray*}
	By using the reproduction property and the Cauchy-Schwarz inequality, we have for all $p\geq 1$,
	$$\mid e_{p}(x)\mid =\mid<e_{p},K(x,\cdot)>_{\mathcal{H}}\mid\leq \parallel e_{p}\parallel_{\mathcal{H}}\,\,\parallel K(x,\cdot)\parallel_{\mathcal{H}}\leq \Vert  K\Vert^{1/2}_{\infty}.$$	
	So we have
\begin{eqnarray*}
\mid Cov(Y^2_{n,p,i,j,l},Y^2_{n,r,i,j,l})\mid &\leq& \mathbb{E}(Y^2_{n,p,i,j,l}Y^2_{n,r,i,j,l}) +\mathbb{E}(Y^2_{n,p,i,j,l})\mathbb{E}(Y^2_{n,r,i,j,l})\nonumber\\&\leq& \mathbb{E}^{1/2}(Y^4_{n,p,i,j,l})\mathbb{E}^{1/2}(Y^4_{n,r,i,j,l})+\mathbb{E}^{1/2}(Y^4_{n,p,i,j,l})\mathbb{E}^{1/2}(Y^4_{n,r,i,j,l}) \nonumber\\&\leq &2 \mathbb{E}^{1/2}(Y^4_{n,p,i,j,l})\mathbb{E}^{1/2}(Y^4_{n,r,i,j,l})\nonumber\\&\leq& C\left[n^{-1}\Vert K\Vert_{\infty} \,\,\mathbb{E}(Y^2_{n,p,i,j,l})\right]^{1/2}\left[n^{-1}\Vert K\Vert_{\infty} \,\,\mathbb{E}(Y^2_{n,r,i,j,l})\right]^{1/2}\nonumber\\&\leq& C n^{-2}\Vert K\Vert_{\infty}\,\, \lambda^{1/2}_{p}\lambda^{1/2}_{r}.
	\end{eqnarray*}
\hfill $\Box$
\subsection{Proof of Theorem \ref{loi}}
From Lemma \ref{l3} it is seen that it suffices to get the asymtotic distribution of $\widetilde{S}_n$.
	We have:
	\begin{eqnarray}\label{som}
	\sum_{j=1}^{k}n_j\parallel V_n^{-1/2}\widehat{\delta}_j \parallel_{\mathcal{H}}^2
&=&\sum_{j=1}^{k}n_j<\widehat{\delta}_j,  V_n^{-1}\widehat{\delta}_j>_{\mathcal{H}}\nonumber\\
&=&\sum_{j=1}^{k}n_j\textrm{tr}\left[\widehat{\delta}_j \otimes( V_n^{-1}\widehat{\delta}_j)\right]\nonumber\\
&=&\sum_{j=1}^{k}n_j\textrm{tr}\left[V_n^{-1}(\widehat{\delta}_j \otimes \widehat{\delta}_j)\right]\nonumber\\
&=&\sum_{j=1}^{k}\sum_{p=1}^{+\infty}n_j<V_n^{-1}e_p,(\widehat{\delta}_j \otimes \widehat{\delta}_j)e_p>_{\mathcal{H}}\nonumber\\
&=&\sum_{j=1}^{k}\sum_{p=1}^{+\infty}(\lambda_p+\gamma_n)^{-1}(\sqrt{n_j}<e_p,\widehat{\delta}_j>_{\mathcal{H}})^2 .
	\end{eqnarray}
	 
	\noindent Considering
	
	\begin{eqnarray*}
	S_{n,p,j}&=&\sqrt{n_j}<\widehat{\delta_j}, e_p>_{\mathcal{H}}\nonumber\\
&=&<\sum_{\stackrel{l=1}{l\neq j}}^{k}P_l\sqrt{n_j}(\widehat{m}_j-\widehat{m}_l),e_p>_{\mathcal{H}} \nonumber\\
&=&\sum_{\stackrel{l=1}{l\neq j}}^{k}P_l\sqrt{n_j}\bigg(\left(<\widehat{m}_j,e_p>_{\mathcal{H}}-<m_j,e_p>_{\mathcal{H}}\right)-\left(<\widehat{m}_l,e_p>_{\mathcal{H}}-<m_l,e_p>_{\mathcal{H}}\right)        \bigg)\nonumber\\
&=&\sum_{\stackrel{l=1}{l\neq j}}^{k}P_l\sqrt{n_j}\left\{\frac{1}{n_j}\sum_{i=1}^{n_j}\left[e_p(X_i^{(j)})-\mathbb{E}[e_p(X_1^{(j)})]\right]-\frac{1}{n_l}\sum_{i=1}^{n_l}\left[e_p(X_i^{(l)})-\mathbb{E}[e_p(X_1^{(l)})]\right]\right\}\nonumber\\
&=&\sum_{\stackrel{l=1}{l\neq j}}^{k}\sum_{i=1}^{n_j}\frac{P_l\sqrt{n_j}}{n_j}\left[e_p(X_i^{(j)})-\mathbb{E}[e_p(X_1^{(j)})]\right]-\sum_{\stackrel{l=1}{l\neq j}}^{k}\sum_{i=n_j+1}^{n_j+n_l}\frac{P_l\sqrt{n_j}}{n_l}\left[e_p(X_{i-n_j}^{(l)})-\mathbb{E}[e_p(X_1^{(l)})]\right]\\
&=&\sum_{\stackrel{l=1}{l\neq j}}^{k}\sum_{i=1}^{n_j+n_l}Y_{n,p,i,j,l} ,  
	\end{eqnarray*}				
we have
\begin{eqnarray*}
S_{n,p,j}^2&=&\sum_{\stackrel{l=1}{l\neq j}}^{k}\sum_{i=1}^{n_j+n_l}Y_{n,p,i,j,l}^2+2\sum_{\stackrel{l=1}{l\neq j}}^{k}\sum_{i=1}^{n_j+n_l}Y_{n,p,i,j,l}\left\{\sum_{t=1}^{i-1}Y_{n,p,t,j,l}\right\}\\
& &+\sum_{\stackrel{l=1}{l\neq j}}^{k}\sum_{\stackrel{q=1}{q\neq l,\,q\neq j}}^{k}\sum_{i=1}^{n_j+n_l}\sum_{i_1=1}^{n_j+n_q}Y_{n,p,i,j,l}Y_{n,p,i_1,j,q}
\end{eqnarray*}
and, from (\ref{som}), 
 \[
\sum_{j=1}^{k}n_j\parallel V^{-1/2}_n\widehat{\delta}_j\parallel^2_{\mathcal{H}}=\sum_{j=1}^k\sum_{p=1}^k(\lambda_p+\gamma_n)^{-1}S_{n,p,j}^{2}=A_n+2B_n+C_n,
\]
where 
	\begin{eqnarray*}
	&&A_{n}=\sum\limits_{p=1}^{+\infty}(\lambda_{p}+\gamma_{n})^{-1}\sum\limits_{j=1}^{k}\sum_{\stackrel{l=1}{l\neq j}}^{k}\sum_{i=1}^{n_j+n_l}\left\{Y^2_{n,p,i,j,l}\right\},\\
	&& B_{n}=\sum\limits_{p=1}^{+\infty}(\lambda_{p}+\gamma_{n})^{-1}\sum\limits_{j=1}^{k}\left\{\sum_{\stackrel{l=1}{l\neq j}}^{k}\sum_{i=1}^{n_j+n_l}Y_{n,p,i,j,l}\left\{\sum\limits_{t=1}^{i-1}Y_{n,p,t,j,l}\right\}\right\},\\
	&& C_n=\sum\limits_{p=1}^{+\infty}\sum\limits_{j=1}^{k}\sum_{\stackrel{l=1}{l\neq j}}^{k}\sum_{\stackrel{q=1}{q\neq l,\,q\neq j}}^{k}\sum_{i=1}^{n_j+n_l}\sum_{i_1=1}^{n_j+n_q}(\lambda_{p}+\gamma_{n})^{-1}Y_{n,p,i,j,l}Y_{n,p,i_1,j,q}.
	\end{eqnarray*}
Thus 
\begin{eqnarray*}
\widetilde{S}_n=2^{-1/2}\ell_n^{-1}\left(A_n+2B_n+C_n\right)
\end{eqnarray*}
and the required result is obtained if we show that  $\ell_n^{-1}A_n$  and $\ell_n^{-1}C_n$  converge in probability to $0$ and that $\ell_n^{-1}B_n$ converges in distribution to $\mathcal{N}(0,1/2)$, as $n\rightarrow +\infty$. These two first properties are obtained if  $A_n=o_P(1)$ and $C_n=o_P(1)$ since $\lim_{n\rightarrow +\infty}\ell_n=+\infty$.

\bigskip

\noindent\textit{Step 1}: let us prove that $A_n=o_P(1)$; it suffices to prove that each term  $B_{j,l}$ defined as $B_{j,l}:=\sum\limits_{p=1}^{+\infty}(\lambda_{p}+\gamma_{n})^{-1}\sum_{i=1}^{n_j+n_l}\left\{Y^2_{n,p,i,j,l}\right\}$ equals $o_P(1)$.  Since  
	$Y_{n,p,i,j,l}$ and $Y_{n,r,t,j,l}$ are independent if $i\neq t$, then $Var(B_{j,k})=\sum\limits_{i=1}^{n_j+n_l}v_{n,i,j,l}$, where 
	\begin{eqnarray}
	v_{n,i,j,l}&=&Var\left( \sum\limits_{p=1}^{+\infty}(\lambda_{p}+\gamma_{n})^{-1}Y^2_{n,p,i,j,l}\right) \nonumber\\
          &=&\sum\limits_{p=1}^{+\infty}\sum\limits_{r=1}^{+\infty}(\lambda_{p}+\gamma_{n})^{-1}(\lambda_{r}+\gamma_{n})^{-1} C_{ov}(Y^2_{n,p,i,j,l},Y^2_{n,r,i,j,l}).\nonumber
	\end{eqnarray} 
	Using equation (\ref{f}), we get 
\[
Var(B_{j,l})\leq C n^{-1}\left(\sum\limits_{p=1}^{+\infty}(\lambda_{p}+\gamma_{n})^{-1}\lambda^{1/2}_{p}\right)^2\leq C n^{-1}\gamma^{-2}_{n}\left(\sum\limits_{p=1}^{+\infty}\lambda^{1/2}_{p}\right)^2
\]
 and since  $\lim_{n\rightarrow +\infty}(\gamma_n^{-1}n^{-1/2})= 0$ and $\sum_{p=1}^{+\infty}\lambda_p^{1/2}<+ \infty$ (because under $\mathscr{H}_0$, $W=V_j$), we   deduce that $B_{j,l}=o_P(1)$. 

\bigskip

\noindent\textit{Step 2}:  let us prove that  $C_n=o_P(1)$. It is easy to check that $C_n=C_{1,n}+C_{2,n}+C_{3,n}+C_{4,n}+C_{5,n}$ where
	\begin{eqnarray}
	&& C_{1,n}=\sum\limits_{j=1}^{k}\sum_{\stackrel{l=1}{l\neq j}}^{k}\sum_{\stackrel{q=1}{q\neq l,\,q\neq j}}^{k}\left[\sum\limits_{p=1}^{+\infty}\sum_{i=1}^{n_j}(\lambda_{p}+\gamma_{n})^{-1}Y_{n,p,i,j,l}Y_{n,p,i,j,q}\right],\nonumber\\
	&& C_{2,n}=\sum\limits_{j=1}^{k}\sum_{\stackrel{l=1}{l\neq j}}^{k}\sum_{\stackrel{q=1}{q\neq l,\,q\neq j}}^{k}\left[\sum\limits_{p=1}^{+\infty}\sum_{i=1}^{n_j}\sum_{i_1=1,\,i\neq i_1}^{n_j}(\lambda_{p}+\gamma_{n})^{-1}Y_{n,p,i,j,l}Y_{n,p,i_1,j,q}\right],\nonumber\\
	&& C_{3,n}= \sum\limits_{j=1}^{k}\sum_{\stackrel{l=1}{l\neq j}}^{k}\sum_{\stackrel{q=1}{q\neq l,\,q\neq j}}^{k}\left[\sum\limits_{p=1}^{+\infty}\sum_{i=1}^{n_j}\sum_{i_1=n_j+1}^{n_j+n_q}(\lambda_{p}+\gamma_{n})^{-1}Y_{n,p,i,j,l}Y_{n,p,i_1,j,q}\right],\nonumber\\
	&& C_{4,n}= \sum\limits_{j=1}^{k}\sum_{\stackrel{l=1}{l\neq j}}^{k}\sum_{\stackrel{q=1}{q\neq l,\,q\neq j}}^{k}\left[\sum\limits_{p=1}^{+\infty}\sum_{i=n_j+1}^{n_j+n_l}\sum_{i_1=1}^{n_j}(\lambda_{p}+\gamma_{n})^{-1}Y_{n,p,i,j,l}Y_{n,p,i_1,j,q}\right],\nonumber\\
	&& C_{5,n}= \sum\limits_{j=1}^{k}\sum_{\stackrel{l=1}{l\neq j}}^{k}\sum_{\stackrel{q=1}{q\neq l,\,q\neq j}}^{k}\left[\sum\limits_{p=1}^{+\infty}\sum_{i=n_j+1}^{n_j+n_l}\sum_{i_1=n_j+1}^{n_j+n_q}(\lambda_{p}+\gamma_{n})^{-1}Y_{n,p,i,j,l}Y_{n,p,i_1,j,q}\right].\nonumber
	\end{eqnarray}
For having $C_n=o_P(1)$ it suffices to prove that $C_{u,n}=o_P(1)$ for $u=1,\cdots,5$. From the equality   $Y_{n,p,i_1,j,q}=P_l^{-1} P_q Y_{n,p,i,j,l}$ we have
	\begin{eqnarray*}
	C_{1,n}=\sum\limits_{j=1}^{k}\sum_{\stackrel{l=1}{l\neq j}}^{k}\sum_{\stackrel{q=1}{l\neq l,\,q\neq j}}^{k}\frac{P_q}{P_l}C_{1,n,j,l}
	\end{eqnarray*}
	where
	$C_{1,n,j,l} =\sum\limits_{p=1}^{+\infty}\sum_{i=1}^{n_j}(\lambda_{p}+\gamma_{n})^{-1}Y_{n,p,i,j,l}^2$.
	Thus it is enough to show that $C_{1,n,j,l}=o_P(1)$. 
	Since $Y_{n,p,i,j,l}$ et $Y_{n,r,t,j,l}$ are independent if  $i\neq t$ , it follows
	\begin{eqnarray}\label{varc1}
	Var(C_{1,n,j,l})&=&\sum_{i=1}^{n_j}Var\left[\sum_{p=1}^{+\infty}(\lambda_p+\gamma_n)^{-1}Y_{n,p,i,j,l}^{2}\right]\nonumber\\
&=&\sum_{p=1}^{+\infty}\sum_{r=1}^{+\infty}(\lambda_p+\gamma_n)^{-1}(\lambda_r+\gamma_n)^{-1}\left[\sum_{i=1}^{n_j}C_{ov}\left(Y_{n,p,i,j,l}^{2},Y_{n,r,i,j,l}^{2}\right)\right]\nonumber\\&\leq&\sum_{p=1}^{+\infty}\sum_{r=1}^{+\infty}(\lambda_p+\gamma_n)^{-1}(\lambda_r+\gamma_n)^{-1}(n_jn^{-2}C\lambda_P^{1/2}\lambda_r^{1/2})\nonumber\\&\leq& Cn^{-1}\bigg(\sum_{p=1}^{+\infty}(\lambda_p+\gamma_n)^{-1}\lambda_p^{1/2}\bigg)^2\nonumber\\&\leq& Cn^{-1}\gamma_n^{-2}\left(\sum_{p=1}^{+\infty}\lambda_p^{1/2}\right)^2.
	\end{eqnarray}
	Since $\lim_{n\rightarrow +\infty}(n^{-1}\gamma_n^{-2})=0$ and   $\sum_{p=1}^{+\infty}\lambda_p^{1/2}<+\infty$, we deduce from (\ref{varc1}) that  $C_{1,n,j,l}=o_P(1)$ and, consequently, $G_{1,n}=o_P(1)$. From similar reasoning, using the fact that   $Y_{n,p,i,j,l}$ and  $Y_{n,p,i_1,j,q}$ are independent when $l\neq q$,  we also obtain that $C_{u,n}=o_P(1)$ for  $u=2,3,4,5$. 

\bigskip
		  
	\noindent\textit{Step 3}: let us  show that $\ell^{-1}_nB_n\stackrel{\mathscr{D}}{\rightarrow}\mathcal{N}(0,1/2)$,  as $n\rightarrow +\infty$, by using the central limit theorem  for triangular arrays of
	martingale differences  (see \cite{hall}). For $i=1,2,\cdots,n_j+n_l$, we consider
	\begin{eqnarray}\label{z}
	M_{n,p,i,j,l}=\sum\limits_{t=1}^{i}Y_{n,p,t,j,l},\,\,\,\,\,\,\zeta_{n,i,j,l}=\ell_n^{-1}\sum\limits_{p=1}^{+\infty}(\lambda_p+\gamma_n)^{-1}Y_{n,p,i,j,l}M_{n,p,i-1,j,l}
	\end{eqnarray}
	and 
\[
\mathcal{F}_{n,j,l,i}=\sigma\left(Y_{n,p,t,j,l},p\in\{1,\cdots,n_j+n_l\},t\in\{0,\cdots,i\}\right).
\] 
By construction, $\zeta_{n,i,j,l}$ is a martingale increment, that is $\mathbb{E}[\zeta_{n,i,j,l}\arrowvert \mathcal{F}_{n,j,l,i-1}]=0$, and we have
\[
\ell^{-1}_nB_n=\sum\limits_{j=1}^{k}\sum_{\stackrel{l=1}{l\neq j}}^{k}\sum_{i=1}^{n_j+n_l}\zeta_{n,i,j,l}.
\]
Then, from Theorem 3.2 and Corollary 3.1 in \cite{hall}, we will obtain $\ell^{-1}_nB_n\stackrel{\mathscr{D}}{\rightarrow}\mathcal{N}(0,1/2)$ if we show that
\begin{equation}\label{conv1}
 Z_n^2:=\sum\limits_{j=1}^{k}\sum_{\stackrel{l=1}{l\neq j}}^{k}\sum_{i=1}^{n_j+n_l}\mathbb{E}[\zeta^2_{n,i,j,l}|\mathcal{F}_{n,j,l,i-1}]\stackrel{P}{\rightarrow}\frac{1}{2}
\end{equation}
and
\begin{equation}\label{conv2}
\max_{i,j,l}\left(\vert \zeta_{n,i,j,l}\vert\right)\stackrel{P}{\rightarrow}0\,\,\,\textrm{ and }\,\,\,\mathbb{E}\left(\zeta_{n,i,j,l}^2\right)\,\,\,\textrm{is bounded in  }n.
\end{equation}
 
\bigskip

\noindent Proof of (\ref{conv1}):   
	We have:
	\begin{eqnarray}
	\mathbb{E}[\zeta^2_{n,i,j,l}|\mathcal{F}_{n,j,l,i-1}]&=&\mathbb{E}[(\ell_n^{-1}\sum\limits_{p=1}^{+\infty}(\lambda_p+\gamma_n)^{-1}Y_{n,p,i,j,l}M_{n,p,i-1,j,l})^2|\mathcal{F}_{n,j,l,i-1}]\nonumber\\&=&\ell_n^{-2}\sum\limits_{p,r=1}^{+\infty}(\lambda_p+\gamma_n)^{-1}(\lambda_r+\gamma_n)^{-1}\mathbb{E}\left[Y_{n,r,i,j,l}Y_{n,p,i,j,l}M_{n,r,i-1,j,l}M_{n,p,i-1,j,l}|\mathcal{F}_{n,j,l,i-1}\right]\nonumber.
	\end{eqnarray}
	Since  $M_{n,r,i-1,j,l}M_{n,p,i-1,j,l}$ is $\mathcal{F}_{n,j,l,i-1}$-measurable, it follows
	\begin{eqnarray}
	\mathbb{E}[\zeta^2_{n,i,j,k}|\mathcal{F}_{n,j,l,i-1}]&=&\ell_n^{-2}\sum\limits_{p,r=1}^{+\infty}(\lambda_p+\gamma_n)^{-1}(\lambda_r+\gamma_n)^{-1}M_{n,r,i-1,j,k}M_{n,p,i-1,j,k}\mathbb{E}\left[Y_{n,r,i,j,l}Y_{n,p,i,j,l}|\mathcal{F}_{n,j,l,i-1}\right]\nonumber.
	\end{eqnarray}
	Hence 
\begin{eqnarray}
	\mathbb{E}\left[Y_{n,r,i,j,l}Y_{n,p,i,j,l}|\mathcal{F}_{n,j,l,i-1}\right]&=&\mathbb{E}\left[Y_{n,r,i,j,l}Y_{n,p,i,j,l}|Y_{n,p,0,j,l},\cdots,Y_{n,p,i-1,j,l}\right]\nonumber\\&=&\mathbb{E}\left[Y_{n,r,i,j,l}Y_{n,p,i,j,l}\right]\nonumber
	\end{eqnarray}
	because $Y_{n,r,i,j,l}$ and $Y_{n,p,t,j,l}$ are independent for any $i\neq t$. So
	\begin{eqnarray}
	\mathbb{E}[\zeta^2_{n,i,j,l}|\mathcal{F}_{n,j,l,i-1}]=\ell_n^{-2}\sum\limits_{p,r=1}^{+\infty}(\lambda_p+\gamma_n)^{-1}(\lambda_r+\gamma_n)^{-1}M_{n,r,i-1,j,l}M_{n,p,i-1,j,l}\mathbb{E}\left[Y_{n,r,i,j,l}Y_{n,p,i,j,l}\right]\nonumber.
	\end{eqnarray}
	Therefore 
	\begin{eqnarray*}
	 Z_n^2
&=&\ell_n^{-2}\sum\limits_{j=1}^{k}\sum_{\stackrel{l=1}{l\neq j}}^{k}\bigg\{\sum\limits_{p=1}^{+\infty}(\lambda_{p}+\gamma_{n})^{-2}\sum_{i=1}^{n_j+n_l}M^2_{n,p,i-1,j,l}\mathbb{E}(Y^2_{n,p,i,j,l})\\
& &+\sum\limits_{p\neq r}^{+\infty}(\lambda_{p}+\gamma_{n})^{-1}(\lambda_{r}+\gamma_{n})^{-1}\sum_{i=1}^{n_j+n_l}M_{n,p,i-1,j,l}M_{n,r,i-1,j,l}\mathbb{E}(Y_{n,p,i,j,l}Y_{n,r,i,j,l})\bigg\}\nonumber\\
&=& E_{n}+F_{n}\nonumber
	\end{eqnarray*}
	where	
	\begin{eqnarray}
	&& E_{n}=\ell_n^{-2}\sum\limits_{p=1}^{+\infty}(\lambda_{p}+\gamma_{n})^{-2}\sum\limits_{j=1}^{k}\sum_{\stackrel{l=1}{l\neq j}}^{k}\sum_{i=1}^{n_j+n_l}M^2_{n,p,i-1,j,l}\mathbb{E}(Y^2_{n,p,i,j,l}),\nonumber\\
	&& F_n=\ell_n^{-2}\sum\limits_{p,r=1,p\neq r }^{+\infty}(\lambda_{p}+\gamma_{n})^{-1}(\lambda_{r}+\gamma_{n})^{-1}\nonumber\\
& &\hspace{3cm}\times \sum\limits_{j=1}^{k}\sum_{\stackrel{l=1}{l\neq j}}^{k}\sum_{i=1}^{n_j+n_l}M_{n,p,i-1,j,l}M_{n,r,i-1,j,l}\mathbb{E}(Y_{n,p,i,j,l}Y_{n,r,i,j,l})\label{F}.
	\end{eqnarray}
	First,
	\begin{eqnarray*}
	\mathbb{E}(M^2_{n,p,i,j,l})&=&\mathbb{E}\left[\sum_{t=1}^{i}Y_{n,p,t,j,l}^2\right]+\mathbb{E}\left[\sum_{t=1}^{i}\sum_{t_1=1,t_1\neq t}^{i}Y_{n,p,t,j,l}Y_{n,p,t_1,j,l}\right]=\sum\limits_{t=1}^{i}\mathbb{E}(Y^2_{n,p,t,j,l})
	\end{eqnarray*}
 because $\mathbb{E}\left[\sum_{t=1}^{i}\sum_{t_1=1,t_1\neq t}^{i}Y_{n,p,t,j,l}Y_{n,p,t_1,j,l}\right]=0 $  since $ Y_{n,p,t,j,l}$ and $Y_{n,p,t_1,j,l}$ are independent when  $ t\neq t_1$. Hence
\begin{eqnarray*}
\mathbb{E}(E_n)&=&\ell_n^{-2}\sum\limits_{p=1}^{+\infty}(\lambda_{p}+\gamma_{n})^{-2}\sum\limits_{j=1}^{k}\sum_{\stackrel{l=1}{l\neq j}}^{k}\sum_{i=1}^{n_j+n_l}\mathbb{E}(M^2_{n,p,i-1,j,l})\mathbb{E}(Y^2_{n,p,i,j,l})\nonumber\\
&=&\ell_n^{-2}\sum\limits_{p=1}^{+\infty}(\lambda_{p}+\gamma_{n})^{-2}\sum\limits_{j=1}^{k}\sum_{\stackrel{l=1}{l\neq j}}^{k}\sum_{i=1}^{n_j+n_l}\sum\limits_{t=1}^{i-1}\mathbb{E}(Y^2_{n,p,t,j})\mathbb{E}(Y^2_{n,p,i,j,l})\nonumber\\
&=&\frac{\ell_n^{-2}}{2}\sum\limits_{p=1}^{+\infty}(\lambda_{p}+\gamma_{n})^{-2}\left\{\left[\sum\limits_{j=1}^{k}\sum_{\stackrel{l=1}{l\neq j}}^{k}\sum_{i=1}^{n_j+n_l}\mathbb{E}(Y^2_{n,p,i,j,l})\right]^2-\sum\limits_{j=1}^{k}\sum_{\stackrel{l=1}{l\neq j}}^{k}\sum_{i=1}^{n_j+n_l}\mathbb{E}^2(Y^2_{n,p,i,j,l})\right\}\nonumber\\
& &-\ell_n^{-2}\frac{1}{2}\sum\limits_{p=1}^{+\infty}(\lambda_{p}+\gamma_{n})^{-2}\sum_{j=1}^{k}\sum_{\stackrel{k_1\neq k_2}{k_1\neq j,\,k_2\neq j}}^{k}\left[\sum_{i=1}^{n_j+n_{k_1}}\mathbb{E}(Y^2_{n,p,i,j,k_1})\right]\left[\sum_{i=1}^{n_j+n_{k_2}}\mathbb{E}(Y^2_{n,p,i,j,k_2})\right]\nonumber\\
& &-\ell_n^{-2}\frac{1}{2}\sum\limits_{p=1}^{+\infty}(\lambda_{p}+\gamma_{n})^{-2}\sum_{j_1\neq j_2}^{k}\sum_{\stackrel{l=1}{l\neq j_1}}^{k}\sum_{\stackrel{l=1}{l\neq j_2}}^{k}\left[\sum_{i=1}^{n_{j_1}+n_l}\mathbb{E}(Y^2_{n,p,i,j_1,l})\right]\left[\sum_{i=1}^{n_{j_2}+n_l}\mathbb{E}(Y^2_{n,p,i,j_2,l})\right]
\end{eqnarray*}	
Using Eq.(\ref{d}), we get

\begin{eqnarray}\label{ddd}
&&\frac{1}{2}\sum\limits_{p=1}^{+\infty}(\lambda_{p}+\gamma_{n})^{-2}\left\{\left[\sum\limits_{j=1}^{k}\sum_{\stackrel{l=1}{l\neq j}}^{k}\sum_{i=1}^{n_j+n_l}\mathbb{E}(Y^2_{n,p,i,j,l})\right]^2-\sum\limits_{j=1}^{k}\sum_{\stackrel{l=1}{l\neq j}}^{k}\sum_{i=1}^{n_j+n_l}\mathbb{E}^2(Y^2_{n,p,i,j,l})\right\}\nonumber\\&=&\frac{1}{2}\sum\limits_{p=1}^{+\infty}(\lambda_{p}+\gamma_{n})^{-2}\left\{\left[\left(1+O(n^{-1})\right)\lambda_p\right]^2-\sum\limits_{j=1}^{k}\sum_{\stackrel{l=1}{l\neq j}}^{k}\left(\frac{n_l^4}{n_jn^4}+\frac{n_j^2n_l}{n^4}\right)\lambda_p^2\right\}\nonumber\\&=&\frac{1}{2}\sum\limits_{p=1}^{+\infty}(\lambda_{p}+\gamma_{n})^{-2}\lambda_p^2\left\{\left[1+O(n^{-1})\right]^2+O(n^{-1})\right\}\nonumber\\&=&\frac{1}{2}\sum\limits_{p=1}^{+\infty}(\lambda_{p}+\gamma_{n})^{-2}\lambda_p^2\left\{1+O(n^{-1})\right\}\nonumber\\&=&\frac{1}{2}\ell_n^2[1+O(n^{-1})]	
\end{eqnarray}
Similarly, we have

\begin{eqnarray}\label{dddd}
&&\frac{1}{2}\sum\limits_{p=1}^{+\infty}(\lambda_{p}+\gamma_{n})^{-2}\sum_{j=1}^k\sum_{\stackrel{k_1\neq k_2}{k_1\neq j,k_2\neq j}}^k\left[\sum_{i=1}^{n_j+n_{k_1}}\mathbb{E}(Y^2_{n,p,i,j,k_1})\right]\left[\sum_{i=1}^{n_j+n_{k_2}}\mathbb{E}(Y^2_{n,p,i,j,k_2})\right]\nonumber\\
&=&\frac{1}{2}\sum\limits_{p=1}^{+\infty}(\lambda_{p}+\gamma_{n})^{-2}\sum_{j=1}^{k}\sum_{\stackrel{k_1\neq k_2}{k_1\neq j,k_2\neq j}}^{k}\left(\frac{n_{k_1}(n_{k_1}+n_j)}{n^2}\right)\left(\frac{n_{k_2}(n_{k_2}+n_j)}{n^2}\right)\lambda_p^2\nonumber\\
&=&\frac{1}{2}\sum\limits_{p=1}^{+\infty}(\lambda_{p}+\gamma_{n})^{-2}\lambda_p^2\left[\sum_{j=1}^{k}\sum_{\stackrel{k_1\neq k_2}{k_1\neq j,k_2\neq j}}^{k}\left(P_{k_1}^2P_{k_2}^2+P_{k_2}P_jP_{k_1}^2+P_{k_1}P_jP_{k_2}^2+P_{k_1}P_{k_2}P_j^2\right)\right]\nonumber\\
&=&\frac{1}{2}\sum\limits_{p=1}^{+\infty}(\lambda_{p}+\gamma_{n})^{-2}\lambda_p^2\,O(n^{-1})\hspace*{0.5cm}\textrm{  because } \hspace*{0.3cm}\sum_{j=1}^{k}P_j^2=O(n^{-1})\hspace*{0.3cm} \textrm{ and }\hspace*{0.3cm}\sum_{j=1}^{k}P_j=1 \nonumber\\
&=&\frac{1}{2}\ell_n^2[O(n^{-1})].
\end{eqnarray}
Furthermore
\begin{eqnarray}\label{ddddd}
&&\frac{1}{2}\sum\limits_{p=1}^{+\infty}(\lambda_{p}+\gamma_{n})^{-2}\sum_{j_1\neq j_2}^{k}\sum_{\stackrel{l=1}{l\neq j_1,\,l\neq j_2}}^{k}\left[\sum_{i=1}^{n_{j_1}+n_l}\mathbb{E}(Y^2_{n,p,i,j_1,l})\right]\left[\sum_{i=1}^{n_{j_2}+n_l}\mathbb{E}(Y^2_{n,p,i,j_2,l})\right]\nonumber\\
&=&\frac{1}{2}\sum\limits_{p=1}^{+\infty}(\lambda_{p}+\gamma_{n})^{-2}\sum_{j_1\neq j_2}^{k}\sum_{\stackrel{l=1}{l\neq j_1,\,l\neq j_2}}^{k}\left(\frac{n_l(n_l+n_{j_1})}{n^2}\right)\left(\frac{n_l(n_l+n_{j_2})}{n^2}\right)\lambda_p^2\nonumber\\
&=&\frac{1}{2}\sum\limits_{p=1}^{+\infty}(\lambda_{p}+\gamma_{n})^{-2}\lambda_p^2\,O(n^{-1})\hspace*{0.5cm} \textrm{because } \,\,\,\,\sum_{j=1}^kP_j^2=O(n^{-1})\hspace*{0.3cm}\textrm{ and }\hspace*{0.3cm}\sum_{j=1}^kP_j=1\,\nonumber\\
&=&\frac{1}{2}\ell_n^2[O(n^{-1})].
\end{eqnarray}
From  (\ref{ddd}), (\ref{dddd}), and (\ref{ddddd}), it follows $\mathbb{E}(E_n)=\frac{1}{2}+o(1)$. Now, let prove that  $E_n-\mathbb{E}(E_n)=o_P(1)$. It suffices to show that  $Var(E_n)=o(1)$. Since
\[
E_{n}-\mathbb{E}(E_{n})=\ell_n^{-2}\sum\limits_{p=1}^{+\infty}(\lambda_{p}+\gamma_{n})^{-2}Q_{n,p}
\]
where
\[
Q_{n,p}=\sum\limits_{j=1}^k\sum_{\stackrel{l=1}{l\neq j}}^k\sum_{i=1}^{n_j+n_l}\mathbb{E}[Y^2_{n,p,i,j,l}]\{M^2_{n,p,i-1,j,k}-\mathbb{E}[M^2_{n,p,i-1,j,k}]\}
\]
it follows
\begin{eqnarray}
Var(E_{n})&=&\ell_n^{-4}\sum\limits_{p=1}^{+\infty}\sum\limits_{r=1}^{+\infty}(\lambda_{p}+\gamma_{n})^{-2}(\lambda_{r}+\gamma_{n})^{-2}\mathbb{E}(Q_{n,p}Q_{n,r})\nonumber\\&=&\ell_n^{-4}\sum\limits_{p=1}^{+\infty}(\lambda_{p}+\gamma_{n})^{-4}\mathbb{E}(Q^2_{n,p})+2\ell_n^{-4}\sum\limits_{r=1}^{+\infty}\sum\limits_{p=1}^{r-1}(\lambda_{p}+\gamma_{n})^{-2}(\lambda_{r}+\gamma_{n})^{-2}\mathbb{E}(Q_{n,p}Q_{n,r}).\nonumber
\end{eqnarray}
We know that $\{M^2_{n,p,i,j,l}-\mathbb{E}[M^2_{n,p,i,j,l}]\}_{1\leq i \leq n_j+n_l}$ is a $\mathcal{F}_{n,i}$-adapted martingale (see \cite{harchaoui1} p. 27).
Let $\nu_{n,p,i,j,l}$ denote the  sequence defined by
 $\nu_{n,p,1,j,l}=M^2_{n,p,1,j,l}-\mathbb{E}[M^2_{n,p,1,j,l}]$ and, for all $i\geq 1$,

\begin{eqnarray*}
\nu_{n,p,i,j,l}&=&M^2_{n,p,i,j,l}-\mathbb{E}[M^2_{n,p,i,j,l}]-\{M^2_{n,p,i-1,j,l}-\mathbb{E}[M^2_{n,p,i-1,j,l}]\}\nonumber\\&=&\left[\sum_{t=1}^{i}Y_{n,p,t,j,l}\right]^2-\sum_{t=1}^{i}\mathbb{E}(Y_{n,p,t,j,l}^2)-\left[\sum_{t=1}^{i-1}Y_{n,p,t,j,l}\right]^2+\sum_{t=1}^{i-1}\mathbb{E}(Y_{n,p,t,j,l}^2)\nonumber\\&=&-\mathbb{E}(Y_{n,p,i,j,l}^2)+\left[\sum_{t=1}^{i}Y_{n,p,t,j,l}\right]^2-\left[\sum_{t=1}^{i-1}Y_{n,p,t,j,l}\right]^2\nonumber\\&=&-\mathbb{E}(Y_{n,p,i,j,l}^2)+\left[\sum_{t=1}^{i-1}Y_{n,p,t,j,l}+Y_{n,p,i,j,l}\right]^2-\left[\sum_{t=1}^{i-1}Y_{n,p,t,j,l}\right]^2\nonumber\nonumber\\&=&Y_{n,p,i,j,l}^2-\mathbb{E}(Y_{n,p,i,j,l}^2)+2Y_{n,p,i,j,l}\sum_{t=1}^{i-1}Y_{n,p,t,j,l}\nonumber\\&=&Y^2_{n,p,i,j,l}-\mathbb{E}[Y^2_{n,p,i,j,l}]+2Y_{n,p,i,j,l}M_{n,p,i-1,j,l}.
\end{eqnarray*}
Then 
\begin{eqnarray*}
Q_{n,p}&=&\sum\limits_{j=1}^{k}\sum_{\stackrel{l=1}{l\neq j}}^k\sum_{i=1}^{n_j+n_l}\mathbb{E}[Y^2_{n,p,i,j,l}]\left\{\left(\sum_{t=1}^{i-1}Y_{n,p,t,j,l}\right)^2-\sum_{t=1}^{i-1}\mathbb{E}[Y^2_{n,p,t,j,l}]\right\}\nonumber\\
&=&\sum\limits_{j=1}^k\sum_{\stackrel{l=1}{l\neq j}}^k\sum_{i=1}^{n_j+n_l}\mathbb{E}[Y^2_{n,p,i,j,l}]\left\{\sum_{t=1}^{i-1}Y^2_{n,p,t,j,l}+2\sum_{t=1}^{i-1}Y_{n,p,t,j,l}\left(\sum_{t_1=1}^{t-1}Y_{n,p,t_1,j,k}\right)\right\}\\
& &-\sum\limits_{j=1}^k\sum_{\stackrel{l=1}{l\neq j}}^k\sum_{i=1}^{n_j+n_l}\mathbb{E}[Y^2_{n,p,i,j,l}]\left\{\sum_{t=1}^{i-1}\mathbb{E}[Y^2_{n,p,t,j,l}]\right\}\nonumber\\
&=& \sum\limits_{j=1}^k\sum_{\stackrel{l=1}{l\neq j}}^k\sum_{i=1}^{n_j+n_l}\mathbb{E}[Y^2_{n,p,i,j,l}]\left\{\sum_{t=1}^{i-1}Y^2_{n,p,t,j,l}+2\sum_{t=1}^{i-1}Y_{n,p,t,j,l}M_{n,p,t-1,j,k}\right\}  \\
& &-\sum\limits_{j=1}^k\sum_{\stackrel{l=1}{l\neq j}}^k\sum_{i=1}^{n_j+n_l}\mathbb{E}[Y^2_{n,p,i,j,l}]\left\{\sum_{t=1}^{i-1}\mathbb{E}[Y^2_{n,p,t,j,l}]\right\}\nonumber\\
&=&\sum\limits_{j=1}^k\sum_{\stackrel{l=1}{l\neq j}}^k\sum_{i=1}^{n_j+n_l}\mathbb{E}[Y^2_{n,p,i,j,l}]\left\{\sum_{t=1}^{i-1}\left(Y^2_{n,p,t,j,l}-\mathbb{E}[Y^2_{n,p,t,j,l}]+2Y_{n,p,t,j,l}M_{n,p,t-1,j,k}\right)\right\}\nonumber\\
&=& \sum\limits_{j=1}^k\sum_{\stackrel{l=1}{l\neq j}}^k\sum_{i=1}^{n_j+n_l}\mathbb{E}[Y^2_{n,p,i,j,l}]\left\{\sum_{t=1}^{i-1}\nu_{n,p,t,j,l}\right\}\nonumber\\
&=&\sum\limits_{j=1}^k\sum_{\stackrel{l=1}{l\neq j}}^k\sum\limits_{i=1}^{n_j+n_l-1}\nu_{n,p,i,j,l}\left[\sum\limits_{t=i+1}^{n_j+n_l}\mathbb{E}[Y^2_{n,p,t,j,l}]\right]\nonumber
\end{eqnarray*}
and, for any $1\leq p\leq r\leq n_j+n_l$,
\begin{eqnarray*}
\mid \mathbb{E}[Q_{n,p}Q_{n,r}]\mid\leq \left(\sum\limits_{j=1}^k\sum_{\stackrel{l=1}{l\neq j}}^k\sum_{t=1}^{n_j+n_l}\mathbb{E}[Y^2_{n,p,t,j,l}]\right)\left(\sum\limits_{j=1}^k\sum_{\stackrel{l=1}{l\neq j}}^k\sum_{t=1}^{n_j+n_l}\mathbb{E}[Y^2_{n,r,t,j,l}]\right)\left\vert\sum_{i=1}^{n_j+n_l-1}\mathbb{E}[\nu_{n,p,i,j,l}\nu_{n,r,i,j,l}]\right\vert.
\end{eqnarray*}
Furthermore,
\begin{eqnarray*}
\mathbb{E}[\nu_{n,p,i,j,l}\nu_{n,r,i,j,l}]&=&\mathbb{E}\bigg[\bigg(Y^2_{n,p,i,j,l}-\mathbb{E}[Y^2_{n,p,i,j,l}]+2Y_{n,p,i,j,l}M_{n,p,i-1,j,l}\bigg)\bigg(Y^2_{n,p,i,j,l}-\mathbb{E}[Y^2_{n,r,i,j,l}]\\
& &+2Y_{n,r,i,j,l}M_{n,r,i-1,j,l}\bigg)\bigg]\nonumber\\&=&\mathbb{E}\left[\left(Y^2_{n,p,i,j,l}-\mathbb{E}[Y^2_{n,p,i,j,l}]\right) \left(Y^2_{n,p,i,j,l}-\mathbb{E}[Y^2_{n,p,i,j,l}]\right)\right]\\
& &+4\mathbb{E}\left[Y_{n,p,i,j,l}M_{n,p,i-1,j,l}Y_{n,r,i,j,l}M_{n,r,i-1,j,l}\right]\nonumber\\
& &+2\mathbb{E}\left[\left(Y^2_{n,p,i,j,l}-\mathbb{E}[Y^2_{n,p,i,j,l}]\right)Y_{n,r,i,j,l}M_{n,r,i-1,j,l}\right]\\
& &+2\mathbb{E}\left[Y_{n,p,i,j,l}M_{n,p,i-1,j,k}\left(Y^2_{n,p,i,j,l}-\mathbb{E}[Y^2_{n,r,i,j,l}]\right)\right]\nonumber\\
&=&Cov(Y^2_{n,p,i,j,l},Y^2_{n,r,i,j,l})+4\mathbb{E}\left[Y_{n,p,i,j,l}Y_{n,r,i,j,l}\right]\mathbb{E}\left[M_{n,p,i-1,j,l}M_{n,r,i-1,j,l}\right]\nonumber\\
& &+2\mathbb{E}\left[\left(Y^2_{n,p,i,j,l}-\mathbb{E}[Y^2_{n,p,i,j,l}]\right)Y_{n,r,i,j,l}M_{n,r,i-1,j,l}\right]\\
& &+2\mathbb{E}\left[Y_{n,p,i,j,l}M_{n,p,i-1,j,l}\left(Y^2_{n,p,i,j,l}-\mathbb{E}[Y^2_{n,r,i,j,l}]\right)\right]\nonumber.
\end{eqnarray*}
Since
 \begin{eqnarray*}
&&\mathbb{E}\left[\left(Y^2_{n,p,i,j,l}-\mathbb{E}[Y^2_{n,p,i,j,l}]\right)Y_{n,r,i,j,l}M_{n,r,i-1,j,l}\right]\\
&=&\mathbb{E}\left[Y^2_{n,p,i,j,l}Y_{n,r,i,j,l}M_{n,r,i-1,j,l}\right]-\mathbb{E}[Y^2_{n,p,i,j,l}]\mathbb{E}\left[Y_{n,r,i,j,l}M_{n,r,i-1,j,l}\right]\nonumber\\&=&\sum_{t=1}^{i-1}\mathbb{E}\left[Y^2_{n,p,i,j,l}Y_{n,r,i,j,l}Y_{n,r,t,j,l}\right]-\sum_{t=1}^{i-1}\mathbb{E}[Y^2_{n,p,i,j,l}]\mathbb{E}\left[Y_{n,r,i,j,l}Y_{n,r,t,j,l}\right]\nonumber\\&=&\sum_{t=1}^{i-1}\mathbb{E}\left[Y^2_{n,p,i,j,l}Y_{n,r,i,j,l}\right]\mathbb{E}\left[Y_{n,r,t,j,l}\right]-\sum_{t=1}^{i-1}\mathbb{E}[Y^2_{n,p,i,j,l}]\mathbb{E}\left[Y_{n,r,i,j,l}\right]\mathbb{E}\left[Y_{n,r,t,j,l}\right]\nonumber\\&=&0 
\end{eqnarray*}
and, similarly, $\mathbb{E}\left[Y_{n,p,i,j,l}M_{n,p,i-1,j,k}\left(Y^2_{n,p,i,j,l}-\mathbb{E}[Y^2_{n,r,i,j,l}]\right)\right]=0$, it follows
\begin{eqnarray}\label{nu}
\left\vert\mathbb{E}[\nu_{n,p,i,j,l}\nu_{n,r,i,j,l}]\right\vert&\leq&\left\vert Cov(Y^2_{n,p,i,j,l},Y^2_{n,r,i,j,l})\right\vert\nonumber\\
& &+4\left\vert\mathbb{E}\{Y_{n,p,i,j,l}Y_{n,r,i,j,l}\}\right\vert\,\,\left\vert\mathbb{E}\{M_{n,p,i-1,j,l}M_{n,r,i-1,j,l}\}\right\vert\nonumber\\
&\leq&Cn^{-1}\lambda_p^{1/2}\lambda_r^{1/2}\nonumber\\
& &+4\left\vert\mathbb{E}\{Y_{n,p,i,j,l}Y_{n,r,i,j,l}\}\right\vert\,\,\left\vert\mathbb{E}\{M_{n,p,i-1,j,l}M_{n,r,i-1,j,l}\}\right\vert.
\end{eqnarray}
On the other hand
\begin{eqnarray}\label{emm}
\mathbb{E}\{M_{n,p,i-1,j,l}M_{n,r,i-1,j,l}\}&=&\mathbb{E}\left[\sum_{t=1}^{i-1}\sum_{s=1}^{i-1}Y_{n,p,t,j,l}Y_{n,r,s,j,l}\right]\nonumber\\
&=&\sum_{t=1}^{i-1}\mathbb{E}\left[Y_{n,p,t,j,l}Y_{n,r,t,j,l}\right]+\sum_{t=1}^{i-1}\sum_{s=1,s\neq t}^{i-1}\mathbb{E}\left[Y_{n,p,t,j,l}Y_{n,r,s,j,k}\right]\nonumber\\
&=&\sum_{t=1}^{i-1}\mathbb{E}\left[Y_{n,p,t,j,l}Y_{n,r,t,j,l}\right]+\sum_{t=1}^{i-1}\sum_{s=1,s\neq t}^{i-1}\mathbb{E}\left[Y_{n,p,t,j,l}\right]\mathbb{E}\left[Y_{n,r,s,j,k}\right]\nonumber\\
&=& \sum\limits_{t=1}^{i-1}\mathbb{E}\left[Y_{n,p,t,j,l}Y_{n,r,t,j,l}\right]
\end{eqnarray}
Thus for all  $(j,l) \in \{1,2,\cdots,k\}^2$ with $j\neq l$, using  Lemma \ref{l3},
\begin{eqnarray}\label{ym}
& &\sum_{i=1}^{n_j+n_l} \mathbb{E}[Y_{n,p,i,j,l}Y_{n,r,i,j,l}]\mathbb{E}[M_{n,p,i-1,j,l}M_{n,r,i-1,j,l}]\nonumber\\
&=&\sum_{i=1}^{n_j+n_l} \mathbb{E}[Y_{n,p,i,j,l}Y_{n,r,i,j,l}]\left\{\sum_{t=1}^{i-1}\mathbb{E}[Y_{n,p,t,j,l}Y_{n,r,t,j,l}]\right\}\nonumber\\
&=&\frac{1}{2}\left\{\left(\sum_{i=1}^{n_j+n_l}\mathbb{E}\{Y_{n,p,i,j,l}Y_{n,r,i,j,l}\}\right)^2-\sum_{i=1}^{n_j+n_l}\mathbb{E}^2[Y_{n,p,i,j,l}Y_{n,r,i,j,l}]\right\}  \nonumber\\
&\leq&\frac{1}{2}\left(\sum_{i=1}^{n_j+n_l}\mathbb{E}\{Y_{n,p,i,j,l}Y_{n,r,i,j,l}\}\right)^2\nonumber\\&\leq& C(\lambda_p^{1/2}\lambda_r^{1/2}\delta_{pr})^2
\end{eqnarray}     	
and, from (\ref{nu}) and (\ref{ym}),
\begin{eqnarray*}\label{111}
\mid \mathbb{E}[Q_{n,p}Q_{n,r}]\mid&\leq&\left(\sum\limits_{j=1}^k\sum_{\stackrel{l=1}{l\neq j}}^k\sum_{t=1}^{n_j+n_l}\mathbb{E}[Y^2_{n,p,t,j,l}]\right)\left(\sum\limits_{j=1}^k\sum_{\stackrel{l=1}{l\neq j}}^k\sum_{t=1}^{n_j+n_l}\mathbb{E}[Y^2_{n,r,t,j,l}]\right)\mid\sum_{i=1}^{n_j+n_l-1}\mathbb{E}[\nu_{n,p,i,j,l}\nu_{n,r,i,j,l}]\mid\nonumber\\&\leq&  C \lambda_p\lambda_r\{Cn^{-1}\lambda_p^{1/2}\lambda_r^{1/2}+ C(\lambda_p^{1/2}\lambda_r^{1/2}\delta_{pr})^2\}\nonumber\\&\leq& C\{ n^{-1}\lambda_p^{3/2}\lambda_r^{3/2}+\lambda_p^2\lambda_r^2(\delta_{pr})^2\}.
\end{eqnarray*}
Thus
\begin{eqnarray*}
Var(E_{n})
&\leq& C\ell_n^{-4}\sum\limits_{p=1}^{+\infty}(\lambda_{p}+\gamma_{n})^{-4}\{n^{-1}\lambda_p^{3}
 +\lambda_p^{4} \}\\
& &+C\ell_n^{-4}\sum\limits_{1\leq p<r}^{+\infty}(\lambda_{p}+\gamma_{n})^{-2}(\lambda_{r}+\gamma_{n})^{-2}\{ n^{-1}\lambda_p^{3/2}\lambda_r^{3/2}\}\nonumber\\
&\leq& C\left\{\ell_n^{-2}n^{-1}\gamma_n^{-1}+\ell_n^{-2} \right\}+Cn^{-1}\ell_n^{-4}\left(\sum_{p=1}^{+\infty}\lambda_p\right)^2.
\end{eqnarray*}
Since  $\lim_{n\rightarrow +\infty}(\gamma_n^{-1}n^{-1})=0$, $ \lim\limits_{n\rightarrow +\infty}\ell_n=+\infty $ (see \cite{harchaoui1},  p. 27) and $\sum_{p=1}^{+\infty}\lambda_p<+\infty$ (because under $\mathscr{H}_0$,   $W=V_j$ and $W$ is a Hilbert-Schmidt operator), the preceding inequality implies $\lim\limits_{n\rightarrow +\infty}Var(E_n)=0$ and, therefore,  $E_n-\mathbb{E}(E_n)=o_P(1)$.  It remains to prove that the term $F_n$ defined in (\ref{F}) satisfies $F_n=o_P(1)$. We have
\begin{eqnarray*}
\mathbb{E}(F_{n})&=&\ell_n^{-2}\sum\limits_{p=1 }^{+\infty}\sum_{\stackrel{r=1}{r\neq p}}^{+\infty}(\lambda_{p}+\gamma_{n})^{-1}(\lambda_{r}+\gamma_{n})^{-1}\bigg\{\\
& &\times \sum\limits_{j=1}^k\sum_{\stackrel{l=1}{l\neq j}}^l\left(\sum_{i=1}^{n_l+n_j}\mathbb{E}(M_{n,p,i-1,j,l}M_{n,r,i-1,j,l})\mathbb{E}(Y_{n,p,i,j,l}Y_{n,r,i,j,l})\right)\bigg\}\nonumber
\end{eqnarray*}
and
\begin{eqnarray}
\mathbb{E}(M_{n,p,i-1,j,l}M_{n,r,i-1,j,l})&\leq& \mathbb{E}^{1/2}(M^2_{n,p,i-1,j,l})\mathbb{E}^{1/2}(M^2_{n,r,i-1,j,l})\nonumber.
\end{eqnarray}
Using (\ref{emm}) we have  $\mathbb{E}[M^2_{n,p,i-1,j,k}]=\sum_{t=1}^{i-1}\mathbb{E}[Y^2_{n,p,t,j,l}]$; then (\ref{d}) implies

\begin{eqnarray*}
\mathbb{E}(M_{n,p,i-1,j,l}M_{n,r,i-1,j,l})&\leq&   C \lambda^{1/2}_{p}\lambda^{1/2}_{r}
\end{eqnarray*}
and
\begin{eqnarray*}
\sum\limits_{j=1}^k\sum_{\stackrel{l=1}{l\neq j}}^k\sum_{i=1}^{n_l+n_j}\mathbb{E}[Y_{n,p,i,j,l}Y_{n,r,i,j,l}]\leq C \lambda_p^{1/2}\lambda_r^{1/2}\delta_{pr}.
\end{eqnarray*}
Hence
\begin{eqnarray}
\bigg(\mathbb{E}(\vert  F_{n}\vert ^2)\bigg)^{1/2}&\leq& C\ell_n^{-2}\sum\limits_{p=1 }^{+\infty}\sum_{\stackrel{r=1}{r\neq p}}^{+\infty}(\lambda_{p}+\gamma_{n})^{-1}(\lambda_{r}+\gamma_{n})^{-1}\lambda_p\lambda_r\delta_{pr}\nonumber
\end{eqnarray}\\
and  $\lim\limits_{n\rightarrow +\infty }\mathbb{E}(\vert  F_{n}\vert ^2)=0$.
Consequently,  $F_n=o_P(1)$, and this permits to conclude that    $Z^2_{n}= \frac{1}{2}+ o_{P}(1)$, i.e. $Z_n^2$ converges in probability to $0$   as $n\rightarrow +\infty$.

\bigskip

\noindent Proof of (\ref{conv2}):  It suffices to prove that
$ \mathbb{E}\left(\max\limits_{1\leq i \leq n_j+n_l} \zeta_{n,i,j,l}^2\right)=o(1)$.
Since $\mid Y_{n,p,i,j,l}\mid \leq C n^{-1/2}\Vert K\Vert^{1/2}_{\infty}$ $\mathbb{P}-a.s$, one has for all $j,k \in \{1,2,\cdots,k\}$:
\begin{eqnarray}\label{dd}
\max_{1\leq i \leq n_j+n_l}\mid\zeta_{n,i,j,l}\mid \leq C\ell_n^{-1} n^{-1/2}\sum\limits_{p=1 }^{+\infty}(\lambda_{p}+\gamma_{n})^{-1}\max_{1\leq i \leq n_j+n_l}\mid M_{n,p,i-1,j,l} \mid. 
\end{eqnarray}	
Using Doob's inequality gives $$\mathbb{E}^{1/2}\left(\max_{1\leq i \leq n_j+n_l}\mid M_{n,p,i-1,j,l} \mid^2\right)\leq 2\mathbb{E}^{1/2}[M^2_{n,p,n_j+n_l-1,j,k}]\leq C \lambda^{1/2}_{p}$$ 
 and, from   (\ref{dd}) and  Minkowski's inequality, we have:
 $$\mathbb{E}^{1/2}\left(\max_{1\leq i \leq n_j+n_l}\zeta^2_{n,i,j,l}\right)\leq C \left\{\ell_n^{-1}\gamma^{-1}_{n}n^{-1/2}\sum\limits_{p=1 }^{+\infty}\lambda^{1/2}_{p}\right\}.$$
The properties   $\lim_{n\rightarrow +\infty}(\gamma^{-1}_{n}n^{-1/2})= 0$ and $\sum\limits_{p=1 }^{+\infty}\lambda^{1/2}_{p} < + \infty$ together with the preceding inequality permit ton conclude that  $ \mathbb{E}\left(\max\limits_{1\leq i \leq n_j+n_l} \zeta_{n,i,j,l}^2\right)=o(1)$.


\end{document}